\documentclass[12pt,a4paper]{article}

\usepackage{tensor}

\usepackage{amssymb,amsmath,amsthm,epsfig}
\numberwithin{equation}{section}
\numberwithin{figure}{section}
\numberwithin{table}{section}
\usepackage{mathrsfs}

\usepackage{latexsym, enumerate}
\usepackage{eepic}
\usepackage{epic}
\usepackage{color}
\usepackage{ifpdf}

\usepackage[numbers,sort&compress]{natbib}
\usepackage{natbib}

\usepackage{verbatim}

\usepackage{dsfont}
\usepackage{multirow}
\usepackage{makecell}
\usepackage{bm}
\usepackage{epstopdf}
\usepackage{graphicx}

\usepackage{subfig}

\usepackage{setspace}

\usepackage{tikz}
\usetikzlibrary{patterns}

\graphicspath{{tuii/}{tuii/locatesource/}{tuii/adeNew/}{tuii/testexa/}{tuii/twoDimDatasets/}{tuii/MoGHD/}}

\usepackage{hyperref}
\hypersetup{hypertex=true,
    colorlinks=true,
    linkcolor=blue,
    filecolor=blue,
    urlcolor=blue,
    citecolor=blue,
    bookmarksopen=true,
}
\usepackage{soul}
\usepackage{xcolor}
\usepackage{color}
\sethlcolor{green}
\soulregister\cite7 
\soulregister\citep7 
\soulregister\citet7 
\soulregister\ref7 
\soulregister\pageref7 

\DeclareMathOperator*{\argmin}{arg\,min}

\usepackage[toc,page]{appendix}  

\usepackage{booktabs}
\usepackage{caption}

\usepackage{algorithm}
\usepackage{algpseudocode}

\theoremstyle{plain}
\newtheorem{lem}{Lemma}[section]
\newtheorem{thm}[lem]{Theorem}
\newtheorem{cor}[lem]{Corollary}
\newtheorem{prop}[lem]{Proposition}

\theoremstyle{definition}
\newtheorem{defn}{Definition}[section]

\theoremstyle{remark}
\newtheorem{rem}{Remark}[section]

\newcommand{\diff}{\,\mathrm{d}}

\begin{document}

\title{\LARGE\bf Geometric Optics Approximation Sampling: A Reflector-Induced Transport Map Framework}

\author{
Zejun Sun\thanks
{School of Mathematics, Hunan University, Changsha 410082, China.
Email: sunzejun@hnu.edu.cn}
\and
Guang-Hui Zheng\thanks
{School of Mathematics, Hunan Provincial Key Laboratory of Intelligent Information Processing and Applied Mathematics, Hunan University, Changsha 410082, China.
Email: zhenggh2012@hnu.edu.cn (Corresponding author)}
}

\date{}
\maketitle
\begin{center}{\bf ABSTRACT}
\end{center}\smallskip
In this paper, we propose Geometric Optics Approximation Sampling (GOAS), a reflector-induced transport-map framework for sampling from target measures. Once a reflecting surface is constructed, the associated transport map is explicitly determined by the physical law of reflection. As a concrete realization, we develop a supporting-hyperellipsoid construction that requires only a discrete approximation of the target measure and does not require gradient information of the target density. The formulation accommodates both density-based and sample-based target representations. A softmin smoothing technique is introduced to obtain a smooth approximate transport map from this piecewise hyperellipsoidal construction. We establish well-posedness and stability of the reflector-induced push-forward measure and derive quantitative error estimates in the maximum mean discrepancy metric, and convergence of continuous statistical observables, including fixed-order moments. Numerical experiments on an analytically tractable example, strongly non-Gaussian targets, sample-based target approximations, and Bayesian inverse problems demonstrate the accuracy and flexibility of GOAS.

\vfill


\section{Introduction}\label{sec-intro}

Sampling from complex probability measures is a fundamental task in statistical inference, uncertainty quantification, and Bayesian inverse problems. Markov Chain Monte Carlo (MCMC) methods \citep{robert2004,tierney1994,steve2011,gelman2013} are widely used for this purpose because of their flexibility and broad applicability. However, MCMC methods generate correlated samples, and slow mixing may substantially reduce the effective sample size. This issue becomes particularly significant when the target measure is defined through an expensive forward model, as in many partial differential equations (PDE) constrained Bayesian inverse problems.

Transport-based sampling methods provide an alternative perspective. Instead of constructing a Markov chain, these methods seek a deterministic map that pushes a simple reference measure forward to a prescribed target measure. Once such a map is available, samples from the reference distribution can be transformed into independent samples from the target approximation. This idea underlies measure transport methods \citep{marzouk2016,el2012,parno2022}, normalizing flows \citep{rezende2015,kobyzev2020,papamakarios2021}, and related map-based sampling approaches. A common feature of these methods is that the sampling task is converted into the construction of a transport map. In many existing approaches, this requires prescribing a parametric form of the map, such as a polynomial expansion, triangular map, or neural-network architecture, and then solving an optimization or training problem.
The quality and computational cost of the resulting sampler may therefore depend on the expressiveness of the selected map class and on the numerical procedure used to fit it.

In this paper, we propose \emph{Geometric Optics Approximation Sampling} (GOAS), a reflector-induced transport-map framework for sampling from target measures. The central idea is to construct a reflecting surface whose associated reflection map transports a simple source measure toward a prescribed target measure. Once the reflecting surface has been constructed, the corresponding map is explicitly determined by the physical law of reflection. 
The reflector can therefore be viewed as a geometric representation of a transport map, and the resulting reflection map can be reused to generate push-forward samples from independent draws of the source measure.

The GOAS provides a transport mechanism that is distinct from directly parameterizing the map itself. In GOAS, the map is induced by the geometry of the reflecting surface rather than postulated in advance as a polynomial expansion or a neural-network architecture. The role of geometric optics is thus to supply a physically interpretable map-construction principle: the sampling problem is recast as the construction of a reflector satisfying a prescribed mass-transport condition. At the framework level, GOAS is not tied to one particular numerical solver for this reflector problem; different reflector-construction techniques may, in principle, be incorporated.
The present paper develops one concrete realization of this framework based on supporting hyperellipsoids. We extend the classical supporting-ellipsoid construction for three-dimensional near-field reflector design \citep{kochengin1997,kochengin1998} to a supporting-hyperellipsoid formulation in general dimension. 
This construction is particularly convenient for sampling because it is gradient-free and uses only a discrete approximation of the target measure.

The required target approximation may be obtained in several ways, depending on how the target measure is specified. If an unnormalized target density is available, the approximation can be constructed from density evaluations, quadrature rules, random point sets, low-discrepancy discretizations, or importance-weighted point sets. If target samples or weighted particles are available, they naturally define empirical or particle approximations of the same form. Thus, GOAS can accommodate both density-based and sample-based descriptions of a target measure. This flexibility is useful in Bayesian inference and uncertainty quantification, where posterior information may be available through density evaluations, particles, samples, or preliminary approximations.

Given such a target approximation, the supporting hyperellipsoid method constructs a reflecting surface as the envelope of hyperellipsoids whose foci are the target points. The associated focal parameters are iteratively updated so that the masses induced by the reflection map match the prescribed target weights. This construction provides a concrete realization of the reflector-induced transport-map framework. Since the resulting surface is generally only continuous and may contain nonsmooth junctions between neighboring hyperellipsoidal sheets, we further introduce a softmin smoothing technique. This smoothing step regularizes the reflecting surface and yields a smooth approximate reflection map for sample generation.

The theoretical analysis in this paper is organized around the geometric optics approximation measure, defined as the push-forward of the source measure under the reflector-induced map. We first establish the well-posedness and stability of this measure, which provides the mathematical foundation for using reflection maps as transport maps for sampling framework. We then derive quantitative error estimates in the Maximum Mean Discrepancy (MMD) metric. These estimates decompose the approximation error into contributions from target discretization, numerical evaluation of the reflector-induced masses, reflector construction tolerance, and softmin smoothing.

The main contributions of this work are summarized as follows.
\begin{itemize}

\item
We introduce GOAS as a reflector-induced transport-map framework for sampling from target measures. Once a reflecting surface is constructed, the associated reflection map is explicitly determined by the physical law of reflection and can be reused to generate independent approximate target samples by pushing forward samples from the source measure.

\item
We develop a supporting-hyperellipsoid construction as one GOAS implementation. The method extends the classical supporting-ellipsoid approach, requires no gradient of the target density, and only a discrete approximation of the target measure.
This formulation accommodates both density-based and sample-based target representations, including density evaluations,
quadrature or importance weighting, empirical samples, or particle approximations.

\item
For this implementation, we introduce a softmin smoothing technique to regularize the piecewise hyperellipsoidal reflecting surface and obtain a smooth approximate reflection map by the reflection law. This provides a practical way to generate samples from this smooth approximate push-forward map

\item
We establish the theoretical foundation and error analysis of GOAS at two complementary levels.
At the framework level, we define the geometric optics approximation measure as the push-forward of the source measure under the reflector-induced transport map and study its well-posedness and stability.
At the implementation level, for the supporting-hyperellipsoid construction developed in this paper, we derive quantitative error estimates in the MMD metric.

\item
We demonstrate the proposed framework through numerical experiments on an analytical reflector example, non-Gaussian target measures, a sample-based Gaussian mixture example, and density-based Bayesian inverse problems. These experiments illustrate the construction of reflector-induced transport maps and their use for generating approximate samples from target measures.

\end{itemize}

The remainder of the paper is organized as follows. Section~\ref{sec:geobay} introduces reflector-induced transport maps, defines the geometric optics approximation measure, and presents its well-posedness and stability properties. Section~\ref{sec:algorithm} develops the supporting-hyperellipsoid realization of GOAS, including target-measure discretization, numerical mass evaluation, softmin regularization, and sample generation. Section~\ref{sec:errorEstimation} presents the MMD error analysis and the resulting estimates for statistical observables. Section~\ref{sec:numexa} reports the numerical experiments, and Section~\ref{sec:conclusion} concludes the paper. Additional geometric derivations, proofs, and implementation details are provided in the appendices.

\section{Reflector-induced transport maps}
\label{sec:geobay}

This section formulates the basic transport-map structure underlying GOAS. 
The main object is a reflection map induced by a reflecting surface. 
Once such a surface is constructed, the reflection map pushes a source measure forward to a measure. 
This push-forward measure will be referred to as the geometric optics approximation measure. 
The purpose of this section is to define this measure and state its well-posedness and stability properties for the sampling framework.

\subsection{Source and target measures}
\label{subsec:source_target}

Let
\[
S_+^n:=\{m\in S^n:m_{n+1}>0\}
\]
be the upper unit hemisphere in $\mathbb R^{n+1}$, and let $\varGamma\subset S_+^n$ be a compact source aperture. The target domain $\Omega$ is a subset of the hyperplane
\[
P:=\{z\in\mathbb R^{n+1}:z_{n+1}=h<0\}.
\]
The source measure is
\begin{align}\label{meaI}
\mu_s(\omega)
:=
\int_{\omega} I(m)\,\diff\sigma(m),
\qquad \omega\subset\varGamma,
\end{align}
where $\sigma$ denotes surface measure on $S^n$ and $I\in L^1(\varGamma)$ is the source density. The target measure is
\begin{align}\label{meaL}
\mu_t(w)
:=
\int_w \pi(z)\,\diff\mu(z),
\qquad w\subset\Omega,
\end{align}
where $\mu$ is the $n$-dimensional Lebesgue measure on $P$ and $\pi\in L^1(\Omega)$ is the target density. 


\subsection{Reflection map induced by a reflecting surface}
\label{subsec:reflection_map}

A reflecting surface is represented as a radial graph
\begin{align}\label{RefSur}
    R_{\rho}(m)=m\rho(m),
    \qquad m\in\varGamma,
\end{align}
where $\rho:\varGamma\to(0,\infty)$ is the polar radius. 
For a sufficiently smooth surface, let $\upsilon(m)$ be the unit normal vector at $R_{\rho}(m)$. 
A ray emitted from the source in the direction $m$ is reflected in the direction
\begin{align}\label{RefDir}
    y(m)=m-2(m\cdot \upsilon(m))\upsilon(m),
\end{align}
according to the physical law of reflection. 
The reflected ray intersects the target plane $P$ at
\begin{align}\label{RefMap}
    T_R(m)
    =
    m\rho(m)+y(m)l(m),
    \qquad m\in\varGamma,
\end{align}
where $l(m)$ is the travel distance from the reflecting surface to the target plane along the reflected direction. 
We call $T_R:\varGamma\to\Omega$ the reflection map associated with $R_\rho$.

\begin{figure}[htbp]
  \centering


\tikzset {_hbwj7gujh/.code = {\pgfsetadditionalshadetransform{ \pgftransformshift{\pgfpoint{89.1 bp } { -128.7 bp }  }  \pgftransformscale{1.32 }  }}}
\pgfdeclareradialshading{_1io0q7yh0}{\pgfpoint{-72bp}{104bp}}{rgb(0bp)=(1,1,1);
rgb(0bp)=(1,1,1);
rgb(25bp)=(0.99,0.02,0.02);
rgb(400bp)=(0.99,0.02,0.02)}


\tikzset{
pattern size/.store in=\mcSize,
pattern size = 5pt,
pattern thickness/.store in=\mcThickness,
pattern thickness = 0.3pt,
pattern radius/.store in=\mcRadius,
pattern radius = 1pt}
\makeatletter
\pgfutil@ifundefined{pgf@pattern@name@_dmduwmw05}{
\pgfdeclarepatternformonly[\mcThickness,\mcSize]{_dmduwmw05}
{\pgfqpoint{0pt}{-\mcThickness}}
{\pgfpoint{\mcSize}{\mcSize}}
{\pgfpoint{\mcSize}{\mcSize}}
{
\pgfsetcolor{\tikz@pattern@color}
\pgfsetlinewidth{\mcThickness}
\pgfpathmoveto{\pgfqpoint{0pt}{\mcSize}}
\pgfpathlineto{\pgfpoint{\mcSize+\mcThickness}{-\mcThickness}}
\pgfusepath{stroke}
}}
\makeatother
\tikzset{every picture/.style={line width=0.75pt}} 

\begin{tikzpicture}[x=0.5pt,y=0.5pt,yscale=-1,xscale=1]

\draw (340.94,414.52)
node[
xslant=1.32,
yscale=0.33,
inner sep=0pt
]
{
\includegraphics[width=160pt]{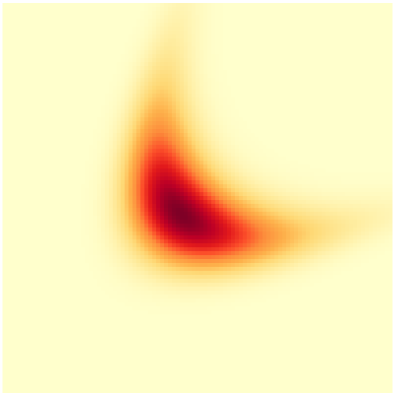}
};
\draw [draw opacity=0][fill={rgb, 255:red, 253; green, 198; blue, 204 }  ,fill opacity=1 ]   (333.58,198.01) .. controls (316.2,178.03) and (271.44,171.98) .. (257.47,173.31) .. controls (243.51,174.63) and (199.08,185.88) .. (188.97,208.42) ;
\draw  [fill={rgb, 255:red, 253; green, 198; blue, 204 }  ,fill opacity=1 ] (260.98,193.01) .. controls (301.1,190.67) and (333.91,193.77) .. (334.26,199.94) .. controls (334.62,206.11) and (302.37,213.01) .. (262.25,215.35) .. controls (222.13,217.69) and (189.32,214.59) .. (188.97,208.42) .. controls (188.62,202.25) and (220.86,195.35) .. (260.98,193.01) -- cycle ;

\draw    (323.41,226.64) -- (146.16,338.11) ;
\draw [shift={(144.47,339.18)}, rotate = 327.83] [color={rgb, 255:red, 0; green, 0; blue, 0 }  ][line width=0.75]    (10.93,-3.29) .. controls (6.95,-1.4) and (3.31,-0.3) .. (0,0) .. controls (3.31,0.3) and (6.95,1.4) .. (10.93,3.29)   ;
\draw    (138.3,263.59) -- (411.13,263.59) ;
\draw [shift={(413.13,263.59)}, rotate = 180] [color={rgb, 255:red, 0; green, 0; blue, 0 }  ][line width=0.75]    (10.93,-3.29) .. controls (6.95,-1.4) and (3.31,-0.3) .. (0,0) .. controls (3.31,0.3) and (6.95,1.4) .. (10.93,3.29)   ;
\draw   (166.63,263.29) .. controls (166.63,213.4) and (210.75,172.96) .. (265.17,172.96) .. controls (319.59,172.96) and (363.71,213.4) .. (363.71,263.29) .. controls (363.71,313.18) and (319.59,353.63) .. (265.17,353.63) .. controls (210.75,353.63) and (166.63,313.18) .. (166.63,263.29) -- cycle ;
\draw  [draw opacity=0] (363.2,265.33) .. controls (362.38,275.84) and (318.55,283.89) .. (264.71,283.35) .. controls (210.81,282.82) and (167.14,273.88) .. (166.65,263.34) -- (264.93,264.09) -- cycle ; \draw   (363.2,265.33) .. controls (362.38,275.84) and (318.55,283.89) .. (264.71,283.35) .. controls (210.81,282.82) and (167.14,273.88) .. (166.65,263.34) ;
\draw  [draw opacity=0][dash pattern={on 4.5pt off 4.5pt}] (166.65,263.24) .. controls (167.43,252.75) and (211.22,244.55) .. (265.06,244.88) .. controls (318.97,245.21) and (362.68,253.96) .. (363.21,264.47) -- (264.93,264.09) -- cycle ; \draw  [dash pattern={on 4.5pt off 4.5pt}] (166.65,263.24) .. controls (167.43,252.75) and (211.22,244.55) .. (265.06,244.88) .. controls (318.97,245.21) and (362.68,253.96) .. (363.21,264.47) ;

\draw [color={rgb, 255:red, 254; green, 142; blue, 142 }  ,draw opacity=1 ]   (342.72,51.35) -- (423.83,413.57) ;
\draw [shift={(424.27,415.52)}, rotate = 257.38] [color={rgb, 255:red, 254; green, 142; blue, 142 }  ,draw opacity=1 ][line width=0.75]    (10.93,-3.29) .. controls (6.95,-1.4) and (3.31,-0.3) .. (0,0) .. controls (3.31,0.3) and (6.95,1.4) .. (10.93,3.29)   ;
\draw  [pattern=_dmduwmw05,pattern size=6pt,pattern thickness=0.75pt,pattern radius=0pt, pattern color={rgb, 255:red, 155; green, 155; blue, 155}][line width=0.75]  (185.76,267.43) .. controls (185.39,261.91) and (217.46,255.3) .. (257.39,252.65) .. controls (297.31,250.01) and (329.97,252.33) .. (330.34,257.85) .. controls (330.7,263.36) and (298.64,269.98) .. (258.71,272.62) .. controls (218.78,275.27) and (186.12,272.94) .. (185.76,267.43) -- cycle ;
\draw  [color={rgb, 255:red, 0; green, 0; blue, 0 }  ,draw opacity=1 ][fill={rgb, 255:red, 246; green, 235; blue, 216 }  ,fill opacity=1 ][line width=0.75]  (220.51,65.29) .. controls (237.97,48.45) and (284.35,19.85) .. (321.77,8.58) .. controls (359.18,-2.7) and (406.63,5.04) .. (420.41,5.92) .. controls (434.18,6.79) and (544.36,50.22) .. (492.86,82.12) .. controls (441.36,114.02) and (172.5,194.66) .. (155.05,168.08) .. controls (137.59,141.49) and (203.06,82.12) .. (220.51,65.29) -- cycle ;
\draw [color={rgb, 255:red, 0; green, 0; blue, 0 }  ,draw opacity=1 ][fill={rgb, 255:red, 246; green, 235; blue, 216 }  ,fill opacity=1 ][line width=0.75]  [dash pattern={on 4.5pt off 4.5pt}]  (155.05,168.08) .. controls (172.5,143.18) and (275.51,89.21) .. (362.8,68.83) .. controls (450.08,48.45) and (484.29,47.81) .. (505.08,59.97) ;

\draw    (342.72,51.35) -- (364.1,145.17) ;
\draw [shift={(364.54,147.12)}, rotate = 257.16] [color={rgb, 255:red, 0; green, 0; blue, 0 }  ][line width=0.75]    (10.93,-3.29) .. controls (6.95,-1.4) and (3.31,-0.3) .. (0,0) .. controls (3.31,0.3) and (6.95,1.4) .. (10.93,3.29)   ;
\draw [color={rgb, 255:red, 249; green, 149; blue, 149 }  ,draw opacity=1 ]   (265.17,263.29) -- (343.14,51.98) ;
\draw [shift={(343.84,50.11)}, rotate = 110.25] [color={rgb, 255:red, 249; green, 149; blue, 149 }  ,draw opacity=1 ][line width=0.75]    (10.93,-3.29) .. controls (6.95,-1.4) and (3.31,-0.3) .. (0,0) .. controls (3.31,0.3) and (6.95,1.4) .. (10.93,3.29)   ;
\draw    (268.27,411.52) -- (265.42,142.76) ;
\draw [shift={(265.39,140.76)}, rotate = 89.39] [color={rgb, 255:red, 0; green, 0; blue, 0 }  ][line width=0.75]    (10.93,-3.29) .. controls (6.95,-1.4) and (3.31,-0.3) .. (0,0) .. controls (3.31,0.3) and (6.95,1.4) .. (10.93,3.29)   ;
\draw [color={rgb, 255:red, 250; green, 177; blue, 45 }  ,draw opacity=1 ]   (342.72,51.35) -- (344.19,11.13) ;
\draw [shift={(344.27,9.13)}, rotate = 92.1] [color={rgb, 255:red, 250; green, 177; blue, 45 }  ,draw opacity=1 ][line width=0.75]    (10.93,-3.29) .. controls (6.95,-1.4) and (3.31,-0.3) .. (0,0) .. controls (3.31,0.3) and (6.95,1.4) .. (10.93,3.29)   ;
\draw    (266.17,263.29) -- (296.59,179.01) ;
\draw [shift={(297.27,177.13)}, rotate = 109.85] [color={rgb, 255:red, 0; green, 0; blue, 0 }  ][line width=0.75]    (10.93,-3.29) .. controls (6.95,-1.4) and (3.31,-0.3) .. (0,0) .. controls (3.31,0.3) and (6.95,1.4) .. (10.93,3.29)   ;
\draw  [draw opacity=0][fill={rgb, 255:red, 254; green, 83; blue, 104 }  ,fill opacity=1 ] (259.03,262.09) .. controls (259.03,258.23) and (262.12,255.09) .. (265.93,255.09) .. controls (269.74,255.09) and (272.83,258.23) .. (272.83,262.09) .. controls (272.83,265.96) and (269.74,269.09) .. (265.93,269.09) .. controls (262.12,269.09) and (259.03,265.96) .. (259.03,262.09) -- cycle ;

\draw (137.54,308.17) node [anchor=north west][inner sep=0.75pt]   [align=left] {$\displaystyle x_{1}$};
\draw (401.67,237.13) node [anchor=north west][inner sep=0.75pt]   [align=left] {$\displaystyle x_{2}$};
\draw (238.92,138.92) node [anchor=north west][inner sep=0.75pt]   [align=left] {$\displaystyle x_{3}$};
\draw (265.96,264.93) node [anchor=north west][inner sep=0.75pt]   [align=left] {$\displaystyle \mathcal{O}$};
\draw (200,430) node [anchor=north west][inner sep=0.75pt]   [align=left] {$\displaystyle \Omega $};
\draw (211.7,198.49) node [anchor=north west][inner sep=0.75pt]   [align=left] {$\displaystyle \varGamma $};
\draw (218.99,256.38) node [anchor=north west][inner sep=0.75pt]   [align=left] {$\displaystyle \Omega _{\varGamma }$};
\draw (281.5,211.7) node [anchor=north west][inner sep=0.75pt]   [align=left] {$\displaystyle m$};
\draw (344.33,41.82) node [anchor=north west][inner sep=0.75pt]   [align=left] {$\displaystyle m\rho ( m)$};
\draw (399.77,406.28) node [anchor=north west][inner sep=0.75pt]   [align=left] {$\displaystyle z$};
\draw (370.86,132.34) node [anchor=north west][inner sep=0.75pt]   [align=left] {$\displaystyle y$};
\draw (278.73,43.18) node [anchor=north west][inner sep=0.75pt]   [align=left] {$\displaystyle R$};
\draw (334.47,436.8) node [anchor=north west][inner sep=0.75pt]   [align=left] {$ $};
\draw (338,435) node [anchor=north west][inner sep=0.75pt]   [align=left] {$\displaystyle \mu _{t}$};
\draw (240,199) node [anchor=north west][inner sep=0.75pt]   [align=left] {$\displaystyle \mu _{s}$};
\draw (271,357) node [anchor=north west][inner sep=0.75pt]   [align=left] {$\displaystyle h$};
\draw (240,320) node [anchor=north west][inner sep=0.75pt]   [align=left] {$\displaystyle S^2$};
\draw (405,298) node [anchor=north west][inner sep=0.75pt]   [align=left] {$\displaystyle l( m)$};
\draw (351,8) node [anchor=north west][inner sep=0.75pt]   [align=left] {$\displaystyle \upsilon $};
\end{tikzpicture}

\caption{
    Schematic illustration of the reflector-induced transport map used in GOAS. 
    A source sample \(m\sim \mu_s\) is mapped to the target domain through the reflection map
    \(z=T_R(m)=m\rho(m)+y(m)l(m)\), where the map is induced by the reflecting surface \(R\) and determined by the physical law of reflection. 
    }
\label{RefFig}
\end{figure}

The key point for GOAS is that once the reflecting surface has been constructed, the map \(T_R\) is explicitly determined by the reflection law. 
Thus the sampling procedure is reduced to drawing samples from the source measure and applying the reflection map.
See Figure \ref{RefFig} for an illustration of sampling a non-Gaussian target distribution with the reflector-induced transport map.

\subsection{Geometric optics approximation measure}
\label{subsec:goa_measure}

We now define the push-forward measure induced by a reflecting surface. 
For nonsmooth or piecewise smooth reflecting surfaces, it is convenient to introduce the visibility sets through supporting hyperellipsoids.

For $z\in\mathbb R^{n+1}$, let $E_z(d)$ denote a hyperellipsoid of revolution with foci $\mathcal O$ and $z$, and focal parameter $d>0$. 
Its polar radius is given by
\begin{align}\label{EllRad}
    \rho_z(x)
    =
    \frac{d}{1-\varepsilon(\hat z\cdot x)},
    \qquad x\in S^n,
\end{align}
where $\hat z=z/|z|$ and
\[
    \varepsilon
    =
    \sqrt{1+\frac{d^2}{|z|^2}}
    -
    \frac{d}{|z|}.
\]
The optical property of such a hyperellipsoid is that rays emitted from one focus are reflected toward the other focus.

\begin{defn}[Supporting hyperellipsoid]\label{def:SE}
Let $R_\rho$ be a reflecting surface. 
A hyperellipsoid $E_z(d)$ with $z\in\Omega$ supports $R_\rho$ at $m\in\varGamma$ if
\[
    \rho(m)=\rho_z(m),
    \qquad
    \rho(x)\leq \rho_z(x),\quad \forall x\in\varGamma .
\]
\end{defn}

Given a reflecting surface $R_\rho$, define the set-valued reflection relation
\[
    \mathcal T_R(x)
    =
    \left\{
    z\in\Omega:
    E_z \text{ supports } R_\rho \text{ at } x\rho(x)
    \right\},
    \qquad x\in\varGamma.
\]
For $z\in\Omega$, define the corresponding visibility set
\[
    V_R(z)
    =
    \left\{
    x\in\varGamma:
    z\in\mathcal T_R(x)
    \right\}.
\]
For a Borel set $w\subset\Omega$, let
\[
    V_R(w)=\bigcup_{z\in w} V_R(z).
\]
At points where $\rho$ is differentiable and the reflection map is single-valued, $\mathcal T_R(x)$ coincides with the reflection map $T_R(x)$ defined in \eqref{RefMap}. 
If $T_R$ is a single-valued and measurable, then $V_R=T_R^{-1}$.

We define the reflector-induced measure on $\Omega$ by
\begin{align}\label{Gomega}
    \mu_R(w)
    :=\mu_s\bigl(V_R(w)\bigr)=
    \int_{V_R(w)} I(x)\,\diff\sigma(x),
    \qquad w\subseteq\Omega .
\end{align}
This measure represents the amount of source mass transported by the reflector to the target set \(w\). 
Equivalently, in the single-valued case,
\[
    \mu_R = (T_R)_{\sharp}\mu_s,
\]
that is $\mu_R(w)=\mu_s(T_R^{-1}(w)),w\subset\Omega$.
Thus, \(\mu_R\) is the push-forward of the source measure under the reflector-induced transport map.

\begin{defn}[Weak solution and geometric optics approximation measure]\label{deGOAM}
A reflecting surface $R_\rho$ is called a weak solution associated with the target measure $\mu_t$ if
\begin{align}\label{GoaMea}
    \mu_R(w)=\mu_t(w),
    \qquad \forall w\subset\Omega
\end{align}
for every Borel set \(w\subset\Omega\). 
The measure \(\mu_R\) defined by \eqref{Gomega} is called the geometric optics approximation measure associated with the reflecting surface \(R_\rho\).
\end{defn}

Equation~\eqref{GoaMea} is the mass-conservation formulation of the near-field reflector problem. For sufficiently smooth reflectors, it is equivalent to the Monge--Amp\`ere-type equation and transport boundary condition summarized in Section~\ref{sec:matref}. From the sampling perspective, Definition~\ref{deGOAM} is the essential formulation: the reflector induces a transport map, and the output distribution is its push-forward of the source measure.
In the present work, we use \eqref{GoaMea} as the measure-theoretic formulation of the reflector-induced transport problem.

\subsection{Well-posedness and stability}
\label{subsec:wellposedness}
The following results provide the theoretical foundation for using reflector-induced maps as transport maps for sampling.

\begin{thm}[Existence and uniqueness of weak solutions]\label{ExiThe}
Let $\Omega\subset P$ be a bounded smooth domain and let $\varGamma\subset S^n$. 
Assume that $I\in L^1(\varGamma)$, $I\geq 0$, and $\pi\in L^1(\Omega)$, $\pi\geq 0$, satisfy the mass condition \eqref{TotEneCon}. 
Let
\[
    \mathcal C_{\varGamma}
    :=
    \left\{
    p\in\mathbb R^{n+1}:
    \frac{p}{|p|}\in\varGamma
    \right\}
\]
be the light cone generated by the source aperture. 
For any point \(p\in\mathcal C_{\varGamma}\) satisfying
\[
    |p|>2\sup_{z\in\Omega}|z|,
\]
there exists a weak solution \(R_\rho\) passing through \(p\). 
Moreover, the solution is unique up to the choice of the normalization point \(p\) and the choice of convex or concave reflector geometry.
\end{thm}

The proof is based on a constructive approximation argument: the target measure is approximated by finite sums of Dirac measures, and the corresponding polyhedral reflectors are constructed as envelopes of hyperellipsoids (see Definition \ref{def:ConRefSur}). 
This classical construction provides the theoretical basis for the supporting hyperellipsoid method developed in Section~\ref{sec:algorithm}; see \citep{kochengin1997,kochengin1998,karakhanyan2010}.

\begin{thm}[Stability of reflecting surfaces]\label{StaRef}
Let \(\{\Omega_k\}_{k\geq1}\) be a sequence of bounded smooth target domains in \(P\) converging to \(\Omega\), and let \(\mu_t^k\) and \(\mu_t\) be the corresponding target measures satisfying the mass condition \eqref{TotEneCon}. 
Assume that \(\mu_t^k\to \mu_t\) as \(k\to\infty\). 
Let \(R_{\rho_k}\) and \(R_\rho\) be the corresponding normalized convex reflecting surfaces. 
Then
\[
    d_H(R_{\rho_k},R_\rho)\to 0,
    \qquad k\to\infty,
\]
where \(d_H(\cdot,\cdot)\) denotes the Hausdorff distance.
\end{thm}

This result states that small perturbations of the target measure or target domain lead to small perturbations of the corresponding reflecting surface. 
It is important for GOAS because the numerical method constructs reflectors from approximations of the target measure.

\begin{thm}[Stability of the geometric optics approximation measure]\label{GOAMstab}
Let the assumptions of Theorem~\ref{StaRef} hold. 
Let \(\mu_{R_k}\) and \(\mu_R\) be the geometric optics approximation measures induced by \(R_{\rho_k}\) and \(R_\rho\), respectively. 
If \(\mu_t^k\to\mu_t\), then
\[
    \mu_{R_k}\to \mu_R
\]
in the weak sense of measures.
\end{thm}

If the original target measure is supported on an unbounded space, we assume that it is transformed to a bounded target domain through a smooth invertible map.
Combining Theorems~\ref{ExiThe}--\ref{GOAMstab}, the geometric optics approximation measure is well defined as a reflector-induced push-forward measure and depends continuously on the target approximation. 
This provides the mathematical foundation for the GOAS framework. 
The next section gives a concrete numerical realization based on a supporting hyperellipsoid construction and softmin smoothing.

\section{GOAS implementation}
\label{sec:algorithm}

The GOAS framework does not prescribe a unique method for constructing the reflecting surface. This section develops the particular realization used in this paper. The target measure is first represented by weighted target points; a piecewise-hyperellipsoidal reflector is then constructed to match these weights; the induced masses are evaluated numerically; and a softmin regularization produces a smooth explicit reflection map. The construction stage is performed offline, after which the resulting map can be reused to generate samples.

\subsection{Discrete approximation of the target measure}
\label{subsec:target_discretization}

Let
$
M:=\mu_s(\varGamma)
$
denote the total source energy. The approximation of the target measure is represented by the nonnegative weighted point measure
\begin{align}\label{eq:DisTarDen}
\mu_t^K
:=
\sum_{i=1}^K w_i\delta_{z_i},
\qquad
z_i\in\Omega,
\quad
w_i\ge0,
\quad
\sum_{i=1}^K w_i=M.
\end{align}
The last condition is the discrete energy-conservation condition required for the supporting-hyperellipsoid construction.
The locations and weights can be constructed from either a density-based or a sample-based description of the target.

If the target density \(\pi\) is evaluable and \(\{z_i,a_i\}_{i=1}^K\) is a quadrature or design rule, we first define
\[
\widetilde w_i:=a_i\pi(z_i)
\]
and rescale the weights by
\begin{align}\label{eq:density_based_weights}
w_i
:=
M\frac{\widetilde w_i}
{\sum_{j=1}^K\widetilde w_j}.
\end{align}
For an importance-weighted random design with proposal density \(q\), one may take
\[
\widetilde w_i=\frac{\pi(z_i)}{q(z_i)}.
\]
If target samples $z_1,\ldots,z_K$ are available, the empirical approximation is obtained by setting $w_i=M/K$.
More generally, weighted particles from importance sampling, sequential Monte Carlo, or another approximation method can
be rescaled to have total mass \(M\) and inserted directly into \eqref{eq:DisTarDen}.

\subsection{Supporting hyperellipsoid construction}
\label{subsec:supporting_hyperellipsoid}

For \(z\in\Omega\) and \(d>0\), let \(E_z(d)\) denote the hyperellipsoid of revolution with foci \(\mathcal O\) and \(z\).
Its radial function is denoted by
\[
f_z(x;d):=\rho_z(x),
\]
and is given by \eqref{EllRad}.
Given target points \(z_1,\ldots,z_K\), Definition \ref{def:ConRefSur} and Proposition \ref{rhoTconvexRef} imply that a piecewise-hyperellipsoidal reflector can be represented as
\begin{align}\label{Rpolyhedron}
\tilde\rho(x;d)
:=
\min_{1\le i\le K}f_{z_i}(x;d_i),
\qquad
R_K(d)
:=
\left\{
x\tilde\rho(x;d):x\in\varGamma
\right\},
\end{align}
where \(d=(d_1,\ldots,d_K)\in(0,\infty)^K\) is the vector of focal parameters.
For each \(i=1,\ldots,K\), define the corresponding visibility cell by
\begin{align}\label{eq:discrete_visibility}
V_i(d)
:=
\left\{
x\in\varGamma:
f_{z_i}(x;d_i)=\tilde\rho(x;d)
\right\}.
\end{align}
The sets \(V_1(d),\ldots,V_K(d)\) form a measurable partition of \(\varGamma\).
Every ray with direction \(x\in V_i(d)\) is redirected toward the
focus \(z_i\). The corresponding induced mass is
\begin{align}\label{eq:induced_mass}
G_i(d)
:=
\mu_s\bigl(V_i(d)\bigr)
=
\int_{V_i(d)}I(x)\,\diff\sigma(x),
\qquad
i=1,\ldots,K.
\end{align}
Therefore, the discrete measure induced by \(R_K(d)\) is
\begin{align}\label{eq:discrete_induced_measure}
\mu_{R_K(d)}
=
\sum_{i=1}^K G_i(d)\,\delta_{z_i}.
\end{align}
Since the visibility cells partition \(\varGamma\), the induced masses automatically satisfy
\begin{align}\label{eq:induced_mass_conservation}
\sum_{i=1}^K G_i(d)
=
\mu_s(\varGamma)
=
M.
\end{align}
The discrete reflector problem is therefore to determine \(d\) such that
\begin{align}\label{eq:discrete_mass_matching}
G_i(d)=w_i,
\qquad
i=1,\ldots,K.
\end{align}
Because both measures have total mass \(M\), only \(K-1\) equations in
\eqref{eq:discrete_mass_matching} are independent.
We measure the mass-matching error by the relative residual
\[
\mathcal E_{\rm ref}(d)
:=
\frac{|G(d)-w|}{M},
\]
where
$
G(d):=(G_1(d),\ldots,G_K(d)),
$
and
$
w:=(w_1,\ldots,w_K).
$
The construction is terminated when
\begin{align} \label{ErrBonDis}
\mathcal E_{\rm ref}(d)\le\epsilon.
\end{align}

To fix the radial normalization of the reflector, we choose the first
focal parameter as a reference. Let
$
Z:=\max_{1\le i\le K}|z_i|,
$
and define
\[
\gamma_i
:=
\max_{x\in\varGamma}x\cdot\hat z_i,
\qquad
\gamma:=\max_{1\le i\le K}\gamma_i<1.
\]
We set
\[
d_1=\alpha Z,
\qquad
\alpha>1,
\]
and restrict the remaining parameters to
\begin{align}\label{eq:focal_parameter_bounds}
c_l d_1\le d_i\le c_r d_1,
\qquad
i=2,\ldots,K,
\end{align}
where
\[
c_l:=\frac{1-\gamma}{2},
\qquad
c_r
:=
\frac{2}{
1-\gamma_1
\left(
\sqrt{1+\dfrac{d_1^2}{Z^2}}
-\dfrac{d_1}{Z}
\right)
}.
\]
These geometric bounds follow from the classical supporting-ellipsoid construction \citep{kochengin1997,kochengin1998}.
The iteration is initialized by
\begin{align}\label{eq:initial_focal_parameters}
d^{(0)}
=
\left(
d_1,c_r d_1,\ldots,c_r d_1
\right).
\end{align}
At this initialization, the nonreference visibility cells have vanishing source masses, so that all source mass is assigned to the reference cell.
For fixed \(d_j\), \(j\neq i\), decreasing \(d_i\) enlarges \(V_i(d)\) and hence increases \(G_i(d)\). This monotonicity forms the basis of the classical supporting-hyperellipsoid iteration \citep{kochengin1998,fournier2010}: the parameters \(d_2,\ldots,d_K\) are decreased successively until the prescribed masses are matched. The procedure is summarized in
Algorithm~\ref{AlgSE}.

\begin{algorithm}[htbp]
\caption{Piecewise-hyperellipsoidal reflector construction}
\label{AlgSE}
\begin{algorithmic}[1]
\Require
Target points and weights
\(\{(z_i,w_i)\}_{i=1}^K\) satisfying
\(\sum_{i=1}^K w_i=M\);
reference parameter \(d_1\);
bounds \(c_l,c_r\);
initial increments \(\Delta_i^{(0)}>0\);
tolerance \(\epsilon>0\);
step tolerance \(\epsilon_{\rm step}\);
maximum number of iterations \(N_{\max}\)
\Ensure
Focal parameter vector \(d^\ast\)
\State Set
$
\tau:={M\epsilon_{\rm ref}}/{\sqrt{K(K-1)}},
$
and initialize using \eqref{eq:initial_focal_parameters}.
\For{\(k=0,\ldots,N_{\max}-1\)}

    \State Compute induced masses
    $
    G(d^{(k)})
    =
    \bigl(G_1(d^{(k)}),\ldots,G_K(d^{(k)})\bigr).
    $

    \If{
    \(\mathcal E_{\rm ref}(d^{(k)})\le\epsilon_{\rm ref}\)
    or
    \(\|\Delta_{2:K}^{(k)}\|_\infty\le\epsilon_{\rm step}\)
    }
        \State \Return \(d^\ast=d^{(k)}\).
    \EndIf

    \State Define the overshooting index set
    $
    J_+^{(k)}
    :=
    \left\{
    i\in\{2,\ldots,K\}:
    G_i(d^{(k)})>w_i+\tau
    \right\}.
    $

    \If{\(J_+^{(k)}\neq\varnothing\)}
        \For{\(i\in J_+^{(k)}\)}
            \State
            \[
            \Delta_i^{(k+1)}
            =
            \frac12\Delta_i^{(k)},
            \qquad
            d_i^{(k+1)}
            =
            d_i^{(k)}+\Delta_i^{(k+1)}.
            \]
        \EndFor
        \State Keep the remaining components unchanged.
        \State \textbf{continue}
    \EndIf

    \State Define
    $
    J_-^{(k)}
    :=
    \left\{
    i\in\{2,\ldots,K\}:
    G_i(d^{(k)})<w_i+\tau
    \right\}.
    $

    \For{\(i\in J_-^{(k)}\)}
        \State
        \[
        \Delta_i^{(k+1)}
        =
        1.25\,\Delta_i^{(k)},
        \qquad
        d_i^{(k+1)}
        =
        \max\left\{
        c_ld_1,\,
        d_i^{(k)}-\Delta_i^{(k+1)}
        \right\}.
        \]
    \EndFor

    \State Keep all remaining components unchanged.

\EndFor
\State \Return the last parameter vector \(d^\ast\).
\end{algorithmic}
\end{algorithm}

\begin{rem}
At each iteration of Algorithm~\ref{AlgSE}, the induced masses \(G_i(d)\), \(i=1,\ldots,K\), must be reevaluated for the current
focal parameters. Their numerical evaluation is described in the next subsection. The iteration is a numerical realization of the constructive supporting-hyperellipsoid procedure underlying Theorem~\ref{ExiThe}. It is gradient-free with respect to the target distribution, since the focal parameters are updated solely from the mass discrepancies \(G_i(d)-w_i\), without evaluating derivatives of the target density.
\end{rem}

\subsection{Numerical evaluation of the induced masses}
\label{subsec:mass_evaluation}

For a given focal-parameter vector \(d\), the induced mass associated
with \(z_i\) is
\[
G_i(d)
=
\int_{\varGamma}
I(x)\mathbf 1_{V_i(d)}(x)\,\diff\sigma(x),
\]
where \(\mathbf 1_{V_i(d)}\) is the indicator function of \(V_i(d)\).
While exact computation through visibility cells (hyperellipsoid intersections) is prohibitively expensive, we use Monte Carlo ray tracing \citep{peter2008, glassner1989} for efficient approximation.
Draw independent samples (i.e., rays) $x_1,\ldots,x_N$ from the distribution $\sigma$ on $S^n$.
For each sample \(x_\ell\), define
\begin{align}\label{IdeSP}
\iota_d(x_\ell)
:=
\operatorname*{arg\,min}_{1\le j\le K}
f_{z_j}(x_\ell;d_j),
\end{align}
The induced mass is then approximated by the Monte Carlo estimator
\begin{align}\label{ApprG}
\widetilde G_i(d)
:=
\frac{\sigma(\varGamma)}{N}
\sum_{\ell=1}^N
I(x_\ell)
\mathbf 1_{\{\iota_d(x_\ell)=i\}},
\qquad
i=1,\ldots,K.
\end{align}

In Algorithm~\ref{AlgSE}, \(G(d)\) is replaced by
\[
\widetilde  G(d)
=
\bigl(\widetilde  G_1(d),\ldots,\widetilde G_K(d)\bigr).
\]
The same source samples are reused throughout the iteration. Each evaluation of \(\widetilde G(d)\) requires \(\mathcal O(KN)\) operations.

\subsection{Softmin smoothing and the approximate reflection map}
\label{subsec:softmin}

The lower envelope \eqref{Rpolyhedron} is continuous but generally not differentiable on the interfaces of the visibility cells. Since the reflection law depends on the surface normal, direct evaluation of the corresponding map is unstable near these interfaces. We therefore propose a softmin smoothing technique \citep{Schmitzer2019}.

\begin{thm}[Softmin smoothing]\label{thm:rhoSmo}
Let $R_K=(d_1,d_2,\dots,d_K)$ be a reflecting surface with a target point set $\{z_i\}_{i=1}^K$.
The polar radius of $R_K$ is
\begin{align}\label{eq:rhodis}
\tilde{\rho}(x)=\min_{1\leq i\leq K}\frac{d_i}{1-\varepsilon_i(\hat{z_i}\cdot x)},\quad x\in\varGamma,
\end{align}
where $\varepsilon_i=\sqrt{1+d_i^2/z_i^2}-d_i/|z_i|$.
Consider the smooth approximation
\begin{align}\label{eq:rhodisSmooth}
\tilde{\rho}_{\lambda}(x)=-\lambda\log\biggl[\sum_{i=1}^{K}\exp(-\frac{f_{z_i}(x;d_i)}{\lambda})\biggr],
\end{align}
where $f_{z_i}(x;d_i)=d_i/(1-\varepsilon_i(\hat{z}_i\cdot x))$ and the parameter $0<\lambda<1$ controls the smoothing.
Then the following approximation error holds:
\begin{align}\label{ineq:rhoappro}
|\tilde{\rho}_{\lambda}-\tilde{\rho}|\leq\lambda\log m+C_1\lambda e^{-\frac{\tau}{\lambda}},
\end{align}
where $m\geq1$ counts the minimal points with $f_{z_{i_k}}(x;d_{i_k})=\tilde{\rho}(x),k\in I^{\ast}:=\{i_1,\dots,i_m\}$, $\tau$ satisfies $\tau_j:=f_{z_j}(x;d_j)-\tilde{\rho}(x)\geq\tau\geq 0$ for any $j\not\in I^{\ast}$, and positive parameter $C_1$ depends on $m$ and $K$.
Moreover, if $m=1$ we have
\begin{align}\label{ineq:Drhoappro}
|D\tilde{\rho}_{\lambda}-D\tilde{\rho}|\leq C_2e^{-\frac{\tau}{\lambda}},
\end{align}
where $C_2>0$ depends on $K$ and $f$.
\end{thm}
We refer to \eqref{eq:rhodisSmooth} as the \emph{softmin smoothing} of \eqref{eq:rhodis}.
When there is a unique minimum, i.e. $m=1$, the approximation error decays exponentially.
As $\lambda\to 0$, $\tilde{\rho}_{\lambda}\to\tilde{\rho}$ with maximal accuracy but minimal smoothness.
Conversely, as $\lambda\to 1$, smoothness increases but approximation quality degrades.
A suitable choice of $\lambda$ balances smoothness and accuracy.

We now provide the explicit approximate expression for the reflection map.

\begin{prop}[Approximate reflection map]\label{prop:refmap}
Let $R_K=(d_1,d_2,\dots,d_K)$ be a reflecting surface associated with a target point set $\{z_i\}_{i=1}^K$.
For the softmin smooth polar radius $\tilde{\rho}_{\lambda}$ in \eqref{eq:rhodisSmooth} of $R_K$, the approximate reflection map $\widetilde{T}_{\lambda}:\varGamma\to\Omega$ is given by
\begin{align}\label{eq:DisT}
z=\widetilde{T}_{\lambda}(x)=\frac{-2\tilde{\rho}_{\lambda}^2D\tilde{\rho}_{\lambda}}{|D\tilde{\rho}_{\lambda}|^2-(\tilde{\rho}_{\lambda}+D\tilde{\rho}_{\lambda}\cdot x)^2}+\biggl(x-\frac{2\tilde{\rho}_{\lambda}D\tilde{\rho}_{\lambda}}{|D\tilde{\rho}_{\lambda}|^2-(\tilde{\rho}_{\lambda}+D\tilde{\rho}_{\lambda}\cdot x)^2}\biggr)\frac{z_{n+1}}{x_{n+1}},
\end{align}
where the gradient $D\tilde{\rho}_{\lambda}$ under the coordinate system \eqref{CooSys} is
$$
D\tilde{\rho}_{\lambda}(x)=\frac{\sum_i\exp(-f(x,z_i)/\lambda)(\varepsilon_if^2(x,z_i)/d_i)(\hat{z}_i-x\hat{z}_{i,n+1}/x_{n+1})}{\sum_i\exp(-f(x,z_i)/\lambda)}.
$$
\end{prop}
Given samples $x$ from the source distribution, the reflection map \eqref{eq:DisT} generates any number of independent approximate samples $z=\widetilde{T}(x)$ from the target distribution.
We now summarize the Geometric Optics Approximation Sampling (GOAS) method proposed in this work as follows:
\begin{description}\label{alg:GOAS}
    \item \vspace{-10pt}\hspace{-5pt}\rule{16.3cm}{0.05em}
    \item [Geometric optics approximation sampling (GOAS):] \vspace{-13pt}
    \item \vspace{-18pt}\hspace{-5pt}\rule{16.3cm}{0.05em}\vspace{-10pt}
    \item[step1] Compute a reflecting surface $R_K=(d_1,d_2,\dots,d_K)$ using Algorithm \ref{AlgSE}. \vspace{-10pt}
    \item[step2] Apply softmin smoothing to obtain a smooth polar radius $\tilde{\rho}_{\lambda}$ of $R_K$, i.e., \eqref{eq:rhodisSmooth}.\vspace{-10pt}
    \item[step3] Generate approximate target samples via the approximate reflection map \eqref{eq:DisT}.
    \item \vspace{-18pt}\hspace{-5pt}\rule{16.3cm}{0.05em}
\end{description}

\section{Error analysis}
\label{sec:errorEstimation}

In this section, we estimate the discrepancy between the approximate
GOAS measure
\[
\mu_{\rm GOAS}
:=
(\widetilde T_{\lambda})_\sharp\mu_s
\]
and the target measure \(\mu_t\). We use the maximum mean discrepancy (MMD) induced by a reproducing kernel Hilbert space \citep{aronszajn1950,berlinet2011,gretton2012}. MMD compares two measures through their kernel mean embeddings and is particularly suitable for the present setting, since it can be evaluated directly from weighted samples or discrete measures without requiring explicit probability densities. This is consistent with the GOAS construction, in which the target measure is first approximated by a weighted sum of Dirac measures. Under the kernel assumptions specified below, MMD distinguishes the relevant measures and provides a convenient framework for decomposing the total error into the target discretization, reflector construction, and smoothing errors.

Let $\mathcal{H}_k$ be a Reproducing Kernel Hilbert Space associated with a positive definite kernel $k: \Omega \times \Omega \to \mathbb{R}$ with inner product $\langle \cdot, \cdot \rangle_{\mathcal{H}_k}$. 
For any functions $f, g \in \mathcal{H}_k$, i.e., $f(\cdot) = \sum_{i} \alpha_i k(z_i, \cdot)$ and $g(\cdot) = \sum_{j}\beta_j k(z_j, \cdot)$, the inner product in $\mathcal{H}_k$ is defined via the kernel as:
\begin{equation}
    \langle f, g \rangle_{\mathcal{H}_k} = \sum_{i} \sum_{j} \alpha_i \beta_j k(z_i, z_j).
\end{equation}
The norm of a function $f \in \mathcal{H}_k$ is then naturally induced by this inner product $\|f\|_{\mathcal{H}_k} = \sqrt{\langle f, f \rangle_{\mathcal{H}_k}}$.
For any $f \in \mathcal{H}_k$, the reproducing property states that $f(x) = \langle f, k(x, \cdot) \rangle_{\mathcal{H}_k}$. 
\begin{defn}[Maximum Mean Discrepancy, MMD]\label{def:MMD}
For two measures $\mu_1, \mu_2$ on $\Omega$, define
\[
d_{k}(\mu_1, \mu_2) = \sup_{f \in \mathcal{H}_k, \|f\|_{\mathcal{H}_k} \leq 1} \left| \int_{\Omega} f \, d\mu_1 - \int_{\Omega} f \, d\mu_2 \right|.
\]
Equivalently, $d_{k}^2(\mu_1, \mu_2) = \|\mu_{1,k} - \mu_{2,k}\|_{\mathcal{H}_k}^2$, where $\mu_{i,k} = \int_{\Omega} k(\cdot, x) \, d\mu_i(x)$ is the kernel mean embedding.
\end{defn}

Throughout this section, we make the following standard assumptions on the kernel:
\begin{enumerate}
\renewcommand{\labelenumi}{\textbf{(A\arabic{enumi})}}
\renewcommand{\theenumi}{\textbf{(A\arabic{enumi})}}
    \item \label{assump:kernel_bounded} The kernel $k$ is bounded, i.e., there exists $\kappa > 0$ such that $\sup_{x \in \Omega} k(x, x) \le \kappa^2$.
    \item \label{assump:kernel_universal} The kernel $k$ is universal, meaning $\mathcal{H}_k$ is dense in $C(\Omega)$ under the supremum norm.
    \item \label{assump:kernel_lipschitz} The kernel $k$ is $L_k$-Lipschitz continuous, that is $|k(x,y) - k(x',y')| \leq L_k(|x-x'| + |y-y'|)$ for all $x,x',y,y' \in \Omega$.
\end{enumerate}

Common examples include the Gaussian kernel and the Mat\'ern kernels, which are universal under appropriate parameter choices.
The assumption \ref{assump:kernel_universal} ensures that convergence in MMD implies weak convergence of measures.
Under assumption \ref{assump:kernel_lipschitz}, any function $f \in \mathcal{H}_k$ is Lipschitz continuous on $\Omega$ \citep{berlinet2011}.

\begin{thm}[Error estimation]\label{thm:errorBound}
Let $\widetilde{T}_{\lambda}:\varGamma\to\Omega$, given by \eqref{eq:DisT}, be the approximate reflection map under $\widetilde{\rho}_{\lambda}$ in \eqref{eq:rhodisSmooth}.
The source measure and target measure are given by $\mu_s$ in \eqref{meaI} and $\mu_t$ in \eqref{meaL}, respectively.
Under assumptions \ref{assump:kernel_bounded}-\ref{assump:kernel_lipschitz}, there exist constants $C_1, C_2 > 0$ such that:
\begin{equation}
    \mathbb{E} [d_{k}(\mu_{\rm GOAS}, \mu_t)] \le \frac{\kappa}{\sqrt{N}} + C_1 (\lambda e^{-1/\lambda} + e^{-1/\lambda}) + C_2 \epsilon  + \mathbb{E} d_k(\mu_t^K,\mu_t).
\end{equation}
where $\epsilon>0$ is any error bound from \eqref{ErrBonDis}.
\end{thm}


\begin{cor}\label{cor:distarget}
If $\mu^K_t$ is constructed by Monte Carlo or importance sampling and the importance weights have finite second moment, then
$$
\mathbb{E} d_k(\mu_t^K,\mu_t)\le C \kappa K^{-1/2},
$$
where $C$ is a positive constant.
\end{cor}

\begin{rem}
The error bound in Theorem \ref{thm:errorBound} quantifies the approximation quality between the geometric optics measure and target measure in Maximum Mean Discrepancy. This bound comprises four distinct components: (1) statistical error from Monte Carlo evaluation of $G(d)$, (2) approximation error due to softmin smoothing ($\lambda$-dependent), (3) prescribed tolerance $\epsilon$ from \eqref{ErrBonDis}, and (4) discretization error of the target measure.
The first and fourth terms vanish as $N,K\to\infty$ (increasing ray samples and target points respectively), while the second and third terms disappear as $\lambda,\epsilon\to 0$ (removing smoothing and tightening tolerances). In the limit, we recover exact agreement between the push-forward of source measure and target measure.
\end{rem}

\begin{thm}\label{thm:moment_error}
Let $\Omega\subset\mathbb R^n$ be compact, and let $\mu_{\rm GOAS}$ and $\mu_t$ be finite positive measures on $\Omega$ with the same total mass $M$. Assume that the kernel $k$ satisfies Assumption~\ref{assump:kernel_universal}. Then, for any $p\in C(\Omega)$ and any $\delta>0$, there exists
$f_\delta\in\mathcal H_k$ such that
\[
\left|
\int_\Omega p(z)\,\diff\mu_{\rm GOAS}(z)
-
\int_\Omega p(z)\,\diff\mu_t(z)
\right|
\le
2M\delta
+
\|f_\delta\|_{\mathcal H_k}
d_k(\mu_{\rm GOAS},\mu_t).
\]
Consequently, if
\[
d_k(\mu_{\rm GOAS},\mu_t)\to0,
\]
then
\[
\int_\Omega p(z)\,\diff\mu_{\rm GOAS}(z)
\longrightarrow
\int_\Omega p(z)\,\diff\mu_t(z)
\]
for every $p\in C(\Omega)$. 
\end{thm}

The above theorem is particularly relevant in applications where one is primarily interested in statistical quantities rather than the full measure itself.


\section{Numerical experiments}\label{sec:numexa}

In this section, we present a series of numerical experiments to illustrate the accuracy, flexibility, and applicability of GOAS. We first consider a spherical reflector for which the exact reflecting surface and target density are analytically available, allowing us to validate the reflector construction and the error behavior with respect to the target-measure discretization and the softmin parameter. We then consider a sample-based target approximation in which the target measure is represented only by a set of weighted samples. Finally, we apply GOAS to two Bayesian inverse problems: acoustic source localization and initial-field reconstruction in a nonlinear advection--diffusion--reaction model. Additional comparisons for strongly non-Gaussian target distributions are provided in Appendix~\ref{subsec:GOASvsMCMCs}.

Unless otherwise stated, all experiments use Algorithm \hyperref[alg:GOAS]{GOAS} with reflector tolerance
\(\epsilon=10^{-4}\) and a uniform source distribution on
\(\varGamma\subset S_+^n\).
Uniform directions on \(S_+^n\) are generated by
\[
X\sim\mathcal N(0,I_{n+1}),\qquad
\widetilde X
=
(X_1,\ldots,X_n,|X_{n+1}|)^\top,
\qquad
x=\frac{\widetilde X}{|\widetilde X|}.
\]

\begin{figure}[htbp]
  \centering
  \subfloat[reflecting surface]
  {   \includegraphics[width=0.11\textwidth,height=0.11\textwidth]{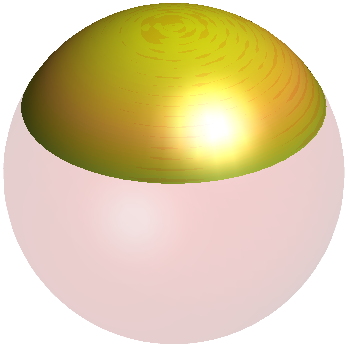}
  }
  \subfloat[K=172]
  {   \includegraphics[width=0.11\textwidth,height=0.11\textwidth]{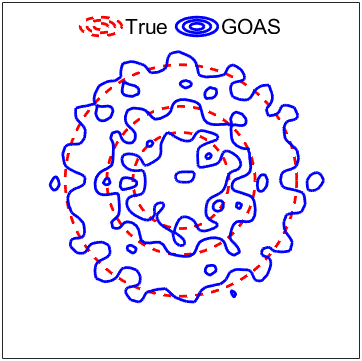}
  }
  \subfloat[K=276]
  {   \includegraphics[width=0.11\textwidth,height=0.11\textwidth]{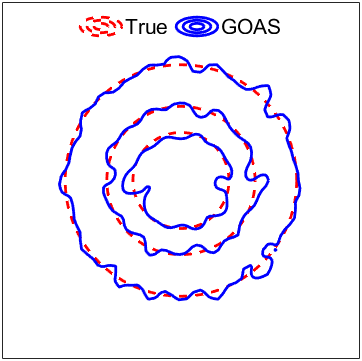}
  }
  \subfloat[K=440]
  {   \includegraphics[width=0.11\textwidth,height=0.11\textwidth]{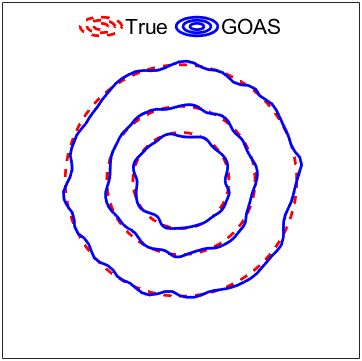}
  }
\subfloat[K=561]
  {   \includegraphics[width=0.11\textwidth,height=0.11\textwidth]{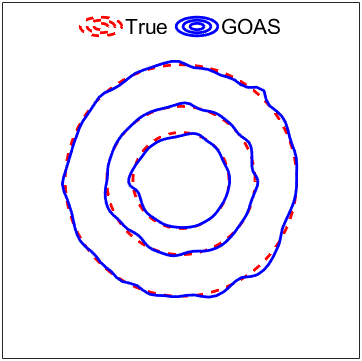}
  }\\
  \subfloat[$\lambda=10^{-4}$]
  {   \includegraphics[width=0.3\textwidth,height=0.2\textwidth]{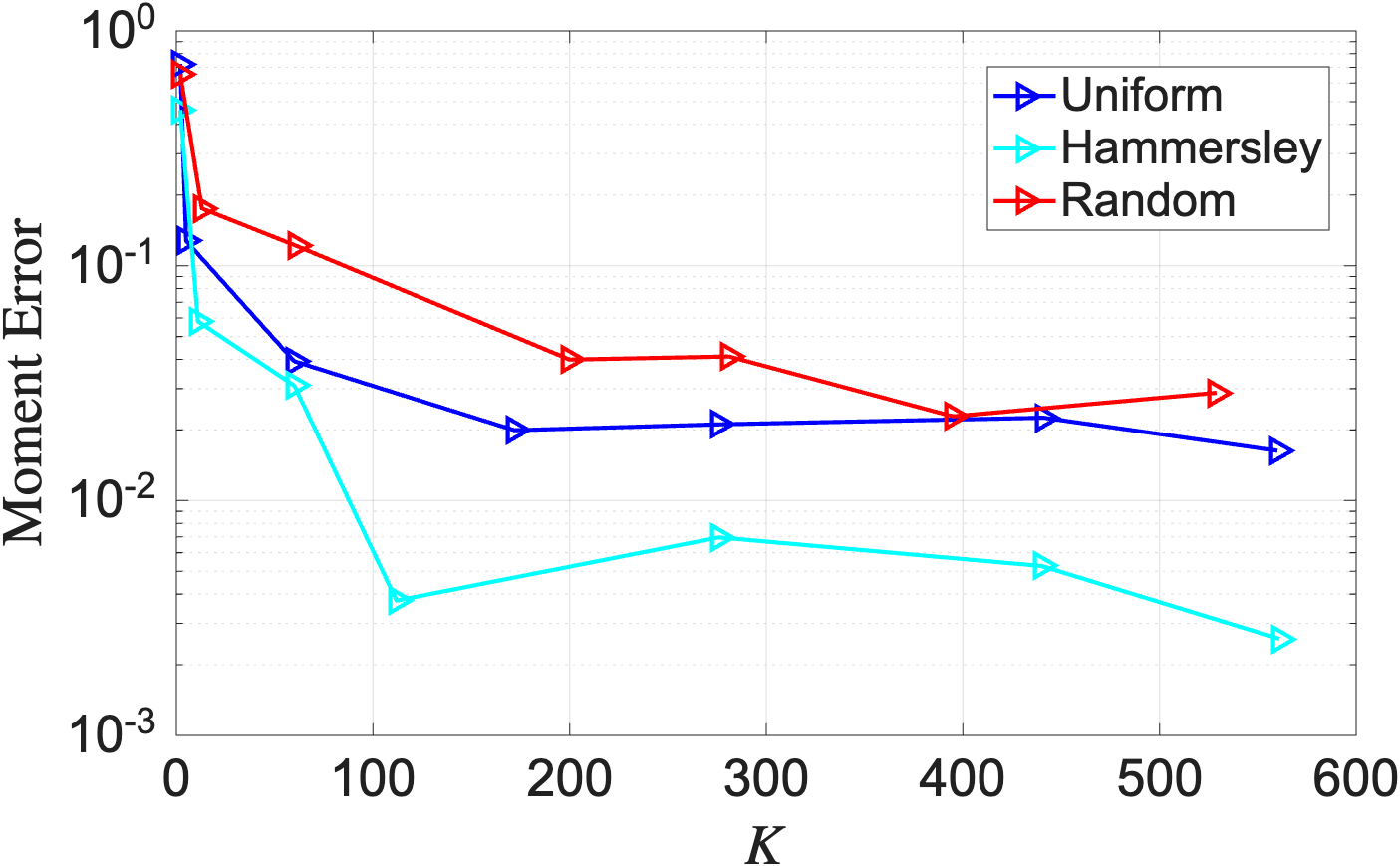}
  }
    \subfloat[K=560]
  {   \includegraphics[width=0.3\textwidth,height=0.2\textwidth]{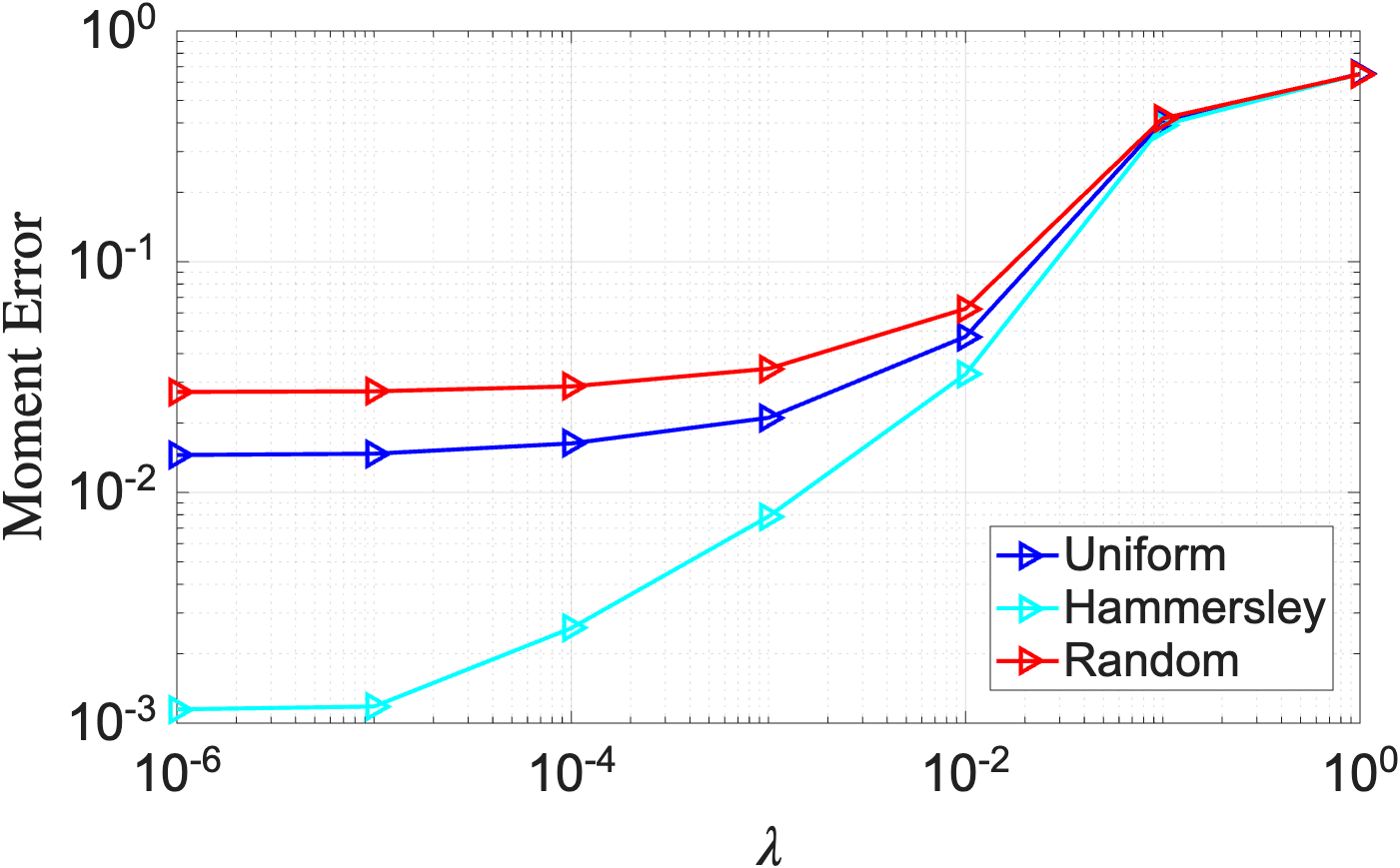}
  }
\caption{
Validation using the spherical reflector.
(a) Reflecting surface reconstructed by the supporting-hyperellipsoid
method with \(K=561\) and \(\lambda=10^{-4}\).
(b)--(e) Exact target density (red) and kernel density estimates of
GOAS samples (blue) for increasing \(K\).
(f) Second-order moment error versus \(K\) for uniform, Hammersley,
and random target point sets with \(\lambda=10^{-4}\).
(g) Second-order moment error versus \(\lambda\) with \(K=560\).
}
 \label{fig:TestGOAS}
\end{figure}
\subsection{Validation using a spherical reflector}\label{sec:TestGOAS}
We present an analytical test case--spherical reflecting surface--to validate our geometric optics approximation sampling method.
Consider $n=2$ with $\Omega\subset P=\{x_3=h|h<0\}$ and $\varGamma\subset S^2_+$.
The reflecting surface is given by
\begin{align}\label{RefSphd}
R(m)=rm,\quad m\in\varGamma,
\end{align}
where $r$ is the polar radius.
The corresponding target density is
\begin{align}\label{LSphd}
\pi(z)=\frac{-hI(-\hat{z})}{|z|^3},\quad z\in\Omega.
\end{align}
When $\varGamma$ is a spherical cap of height $h_c$, $\Omega$ forms a disk with radius $-h\sqrt{{1}/{(1-h_c)^2}-1}$ (see Appendix \ref{AppendixSph} for derivation).

For numerical verification, we set $h=-1,h_c=1-1/\sqrt{5}$ and $r=2$.
Figure \ref{fig:TestGOAS}a shows a spherical sheet with a radius $r=2$, demonstrating the effectiveness of the supporting hyperellipsoid method with softmin smoothing.
The kernel density estimates in Figs.\ref{fig:TestGOAS}b-e show improved approximation to the true density as we increase the number of points in the Hammersley sequence for target distribution discretization.
Figures~\ref{fig:TestGOAS}(f)--(g) further examine the effects of target discretization and softmin smoothing. In Fig.~\ref{fig:TestGOAS}(f), the target measure is discretized using uniform, Hammersley, and random point sets. These point sets provide different discrete approximations of the same target measure, and the second-order moment error generally decreases as \(K\) increases. This numerically illustrates the stability of the reflector construction and the geometric optics approximation measure with respect to target-measure discretization, consistent with Theorems~\ref{StaRef} and~\ref{GOAMstab}.
Figure~\ref{fig:TestGOAS}(g) shows that the error decreases as the softmin parameter \(\lambda\) is reduced, in agreement with the smoothing-error estimate. Together, these results are consistent with the error estimate in Theorem~\ref{thm:errorBound} and the moment
estimate in Theorem~\ref{thm:moment_error}.

\subsection{Sample-based target approximation} \label{subsec:sample_based}

We next consider the case where the target distribution is available
only through samples. The target distribution is chosen as a
\(50\)-dimensional Gaussian mixture
\[
\pi(z)
=
\frac12\mathcal N(z;\mu_1,I_{50})
+
\frac12\mathcal N(z;\mu_2,I_{50}),
\]
where
\[
\mu_1=(0,\ldots,0)^\top,
\qquad
\mu_2=(5,\ldots,5)^\top.
\]
Although the density is explicitly known in this example, it is used
only to generate reference samples and is not used in the construction
of the GOAS reflector.

We first generate \(N_t=10^4\) independent target samples and compress them into \(K=800\) representative points using the \(k\)-medoids algorithm \cite{kaufman2009finding,park2009simple} with squared Euclidean distance. Let \(\mathcal C_i\) denote the cluster associated with the \(i\)-th medoid \(z_i\). We assign to each medoid the empirical cluster weight
\[
w_i
=
\frac{|\mathcal C_i|}{N_t},
\qquad
i=1,\ldots,K,
\]
so that the sample-based target approximation is
\[
\mu_t^K
=
\sum_{i=1}^K w_i\delta_{z_i},
\qquad
\sum_{i=1}^K w_i=1.
\]
Only the medoids \(\{z_i\}_{i=1}^K\) and their weights
\(\{w_i\}_{i=1}^K\) are then supplied to the supporting-hyperellipsoid
construction.

\begin{figure}[htbp]
  \centering
 \includegraphics[width=0.6\textwidth,height=0.4\textwidth]{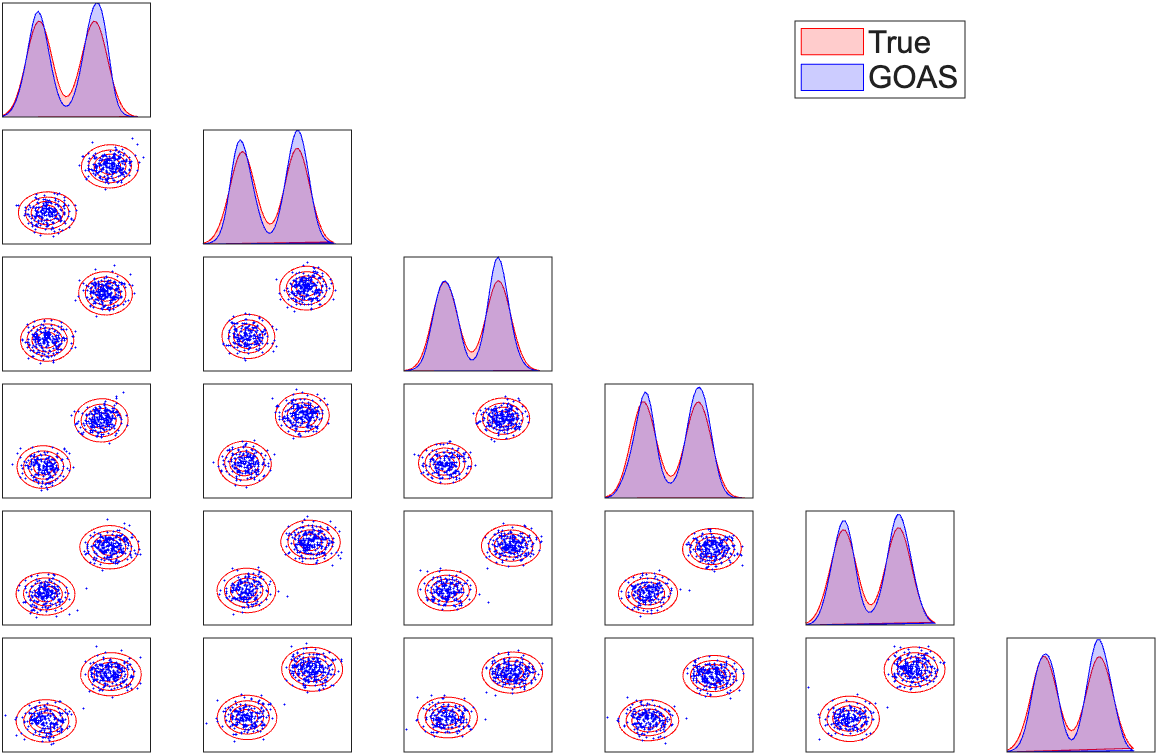}
 \caption{
Comparison of GOAS samples with the target distribution in the sample-based experiment (\(K=800\), \(\lambda=5\times10^{-4}\)). Shown are marginal kernel density estimates and two-dimensional projections for the selected coordinates \(1,10,20,30,40,\) and \(50\).
}
    \label{fig:sample_based}
\end{figure}

Figure~\ref{fig:sample_based} compares the resulting GOAS samples with
independent reference samples from the target distribution. The
one-dimensional marginal kernel density estimates and selected
two-dimensional projections show that GOAS accurately reproduces the
bimodal structure of the target distribution across the selected
coordinates. This example demonstrates that the GOAS construction can
be applied directly to a sample-based representation of the target
measure, without evaluating either the target density or its gradient.

\subsection{Density-based target approximation: Bayesian inverse problems}\label{sec:BayInvPro}

We next consider the density-based setting and demonstrate GOAS on two Bayesian inverse problems: acoustic source localization and initial-field reconstruction in a nonlinear advection--diffusion--reaction model. The corresponding Bayesian formulations are given in Appendix~\ref{Appendix:bayfra}.

In both examples, the posterior density is available pointwise up to a normalizing constant. Since the construction of an optimal discrete approximation of a general posterior measure is not the focus of this work, we employ simple strategies to obtain the discrete target measure
\[
\mu_t^K=\sum_{i=1}^K w_i\delta_{z_i}.
\]
For the acoustic source localization problem, we use an importance-sampling approximation with a Gaussian proposal whose mean and covariance are estimated from a pilot MCMC run. For the nonlinear advection--diffusion--reaction problem, the pilot MCMC samples are used only to identify a bounded region containing the dominant posterior mass; this region is subsequently discretized using a low-discrepancy point set, with weights determined by evaluations of the posterior density. In both cases, neither the normalizing constant nor gradients of the posterior density are required for the subsequent GOAS construction.

\subsubsection{Locating acoustic sources}\label{example1}
Consider the Helmholtz equation
\begin{align}\label{Hel}
\Delta u +k^2u=S \quad\text{in}\ \mathbb{R}^2,
\end{align}
where $k>0$ is the wave number and the field $u$ satisfies the Sommerfeld Radiation Condition, and the source is
$
S(x)=\sum_{i=1}^N\varsigma_i\delta(x-z_i),
$
where $z\in\Omega$ and $\varsigma\neq0$.
The far field pattern is given by
\begin{align}\label{FarFie}
u_{\infty}(\hat{x})=-\frac{e^{i\frac{\pi}{4}}}{\sqrt{8\pi k}}\sum_{i=1}^{N}\varsigma_i e^{-ik\hat{x}\cdot z_i}\quad\hat{x}\in S^1,
\end{align}
where $i=\sqrt{-1}$.
The inverse problem involves recovering the source locations $z_i$ from far field measurements \citep{liu2021,eller2009,el2011}.
More details about the inverse problem can be found in Appendix \ref{Appendix:detailsLAS}.

\begin{figure}[htbp]
  \centering
 \includegraphics[width=0.6\textwidth,height=0.4\textwidth]{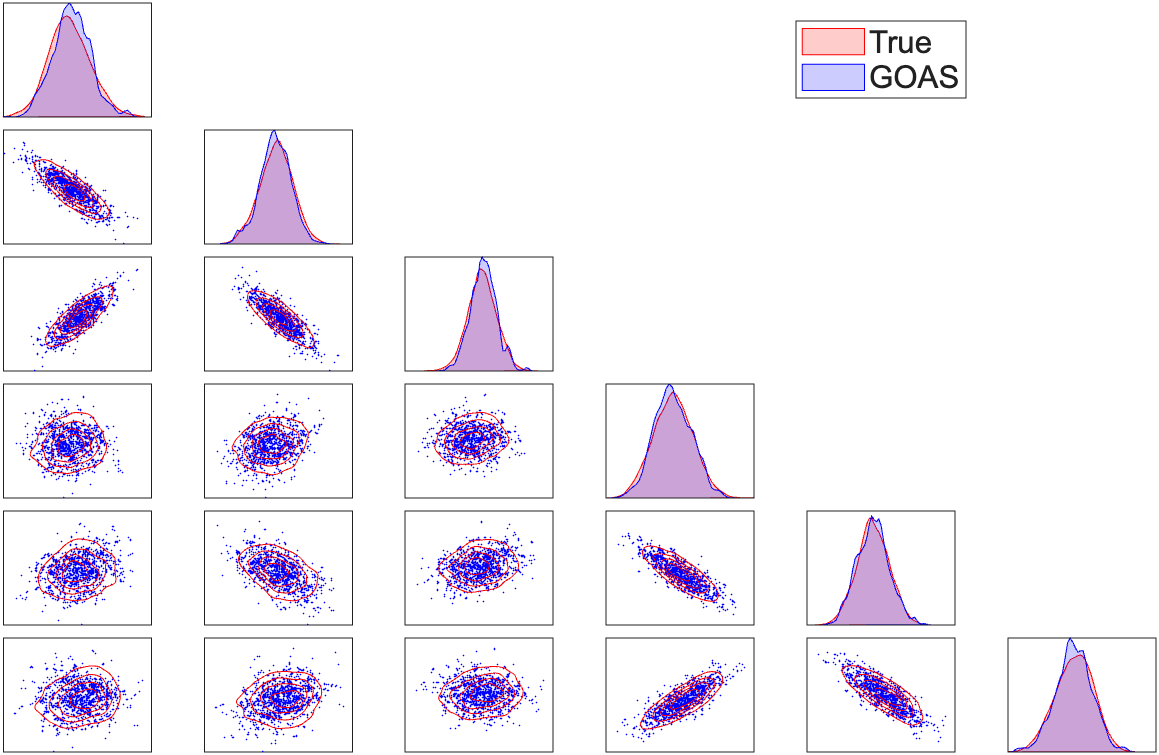}
 \caption{
Comparison of posterior samples obtained from GOAS (\(\lambda=10^{-4}\), \(K=1000\)) and the reference MCMC samples. Diagonal panels show marginal kernel density estimates, while off-diagonal panels show pairwise sample projections.
}
    \label{fig:LAS}
\end{figure}

Table~\ref{tab:LAS} compares the posterior sample means and standard deviations obtained from GOAS with those estimated from the reference MCMC samples. Figure~\ref{fig:LAS} further compares the marginal kernel density estimates and pairwise sample projections. GOAS accurately recovers the posterior means and hence the three source locations, while also capturing the principal marginal and correlation structures of the posterior distribution. Some moderately larger posterior standard deviations are observed for GOAS, but the overall posterior structure remains consistent with the MCMC reference.


\begin{table}[htbp]
\centering
\caption{
Posterior sample means and standard deviations (Std) obtained from GOAS with \(\lambda=10^{-4}\) and \(K=1000\), compared with the reference MCMC estimates.
}
\begin{tabular}{cccccccc}
\hline
\multirow{2}{*}{Locations} &
\multicolumn{2}{c}{MCMC}& &\multicolumn{2}{c}{GOAS} \\
\cline{2-3}\cline{5-6}
& Mean & Std &  & Mean & Std \\ \hline
(1, 4) & (1.041, 4.044) & (0.04150, 0.03917) & & (1.0503, 4.0488) & (0.0426, 0.0411) \\
(2, 5) & (1.994, 4.979) & (0.04521, 0.04288)& &  (1.9916, 4.9763) & (0.0458, 0.0462)  \\
(3, 6)& (3.009, 6.037) & (0.03929, 0.03740) & & (3.0138, 6.0389) & (0.0401, 0.0406) \\ \hline
\end{tabular}
  \label{tab:LAS}
\end{table}

\subsubsection{Initial field reconstructing in nonlinear advection-diffusion-reaction systems}
Consider the advection-diffusion-reaction (ADR) initial-boundary value problem
\begin{align}\label{ADR}
  \frac{\partial u}{\partial t}+\nu\cdot\nabla u-\nabla\cdot(m_d\nabla u)+cu^3&=f\quad \text{in}\;\Omega\times(0,T), \\
  \frac{\partial u}{\partial n} &= 0 \quad\text{on}\;\partial\Omega\times(0,T),\\
  u|_{t=0}&= m_0\quad\text{in}\;\Omega,\label{ADR3}
\end{align}
where $n$ is unit normal of $\partial\Omega$, and $m_d\in L^2(\Omega)$ and $m_0\in L^2(\Omega)$ are the diffusion coefficient and initial field, respectively.
The $\nu$ is advection velocity field, $c$ is the reaction coefficient and $f$ is the source term.
The inverse problem is to reconstruct the initial condition field $m_0$ from terminal-time measurements $u(x,T),x\in\Omega$ \citep{ghattas2021}.

\begin{figure}[htbp]
  \centering
  \subfloat[True]
  {
 \includegraphics[width=0.25\textwidth,height=0.2\textwidth]{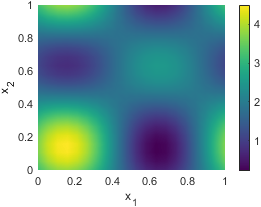}
  }
  \subfloat[MCMC, mean]
  {
    \includegraphics[width=0.26\textwidth,height=0.21\textwidth]{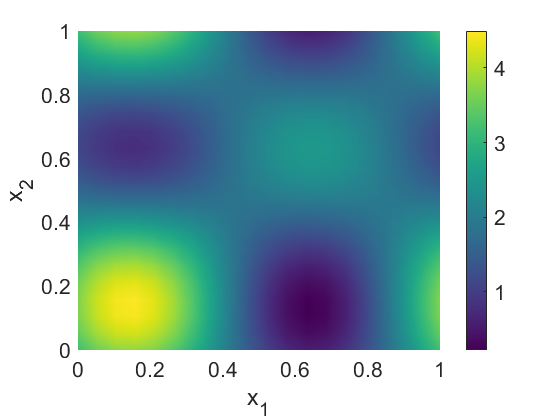}
  }
   \subfloat[MCMC, Std]
  {
  \includegraphics[width=0.26\textwidth,height=0.21\textwidth]{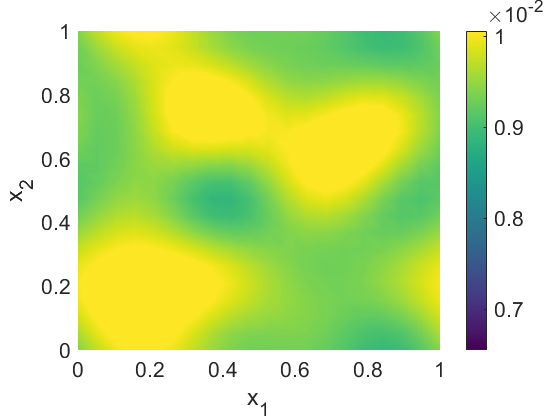}
  }\\
  \subfloat[GOAS, mean]
  {
    \includegraphics[width=0.25\textwidth,height=0.2\textwidth]{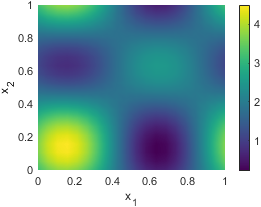}
  }
   \subfloat[GOAS, Std]
  {
  \includegraphics[width=0.26\textwidth,height=0.21\textwidth]{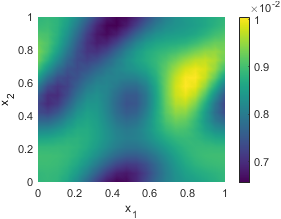}
  }
  \caption{
Comparison of the true initial field with posterior statistics obtained from MCMC and GOAS (\(\lambda=10^{-4}\), \(K=285\)): (a) true initial field; (b)--(c) MCMC posterior mean and standard deviation (Std); (d)--(e) GOAS posterior mean and standard deviation.
}
    \label{fig:ADR}
\end{figure}

Implementation details appear in Appendix \ref{Appendix:detailsADR}.
Using GOAS with \(\lambda=10^{-4}\) and \(K=285\), Figure~\ref{fig:ADR} compares the reconstructed initial field and posterior uncertainty with the MCMC reference. The posterior mean obtained by GOAS closely agrees with the MCMC estimate and accurately recovers the main spatial features of the true initial field. The posterior standard deviations obtained by GOAS and MCMC are also of comparable magnitude, although some differences in their spatial patterns are observed. These results demonstrate the applicability of GOAS to posterior approximation for a nonlinear PDE-constrained inverse problem.

\section{Discussion and conclusions}
\label{sec:conclusion}

\subsection{Discussion}

The GOAS framework is not restricted to the supporting-hyperellipsoid construction developed in this paper. Its essential ingredient is a reflecting surface whose associated transport map is explicitly determined by the physical law of reflection. In principle, other numerical methods for constructing such reflecting surfaces may also be incorporated into GOAS. The following discussion therefore concerns primarily the supporting-hyperellipsoid realization considered in this work, whose numerical properties and limitations need not be intrinsic to the general GOAS framework.

The present implementation requires a discrete approximation
\[
\mu_t^K=\sum_{i=1}^K w_i\delta_{z_i}
\]
of the target measure as input. The construction of this discrete measure is problem dependent and is not the focus of the present work. As illustrated in the numerical experiments, it may be obtained from importance sampling, density-weighted low-discrepancy point sets, or sample-based clustering. Consequently, the overall approximation accuracy depends in part on the quality of this target discretization.

The main computational cost of the present implementation is associated with the repeated evaluation of the reflector-induced masses during the supporting-hyperellipsoid iteration. With \(N\) Monte Carlo samples from the source measure and \(K\) target support points, each mass evaluation requires determining the active supporting hyperellipsoid for every source sample, resulting in a cost of order \(\mathcal O(NK)\). The number \(N\) also controls the Monte Carlo accuracy in the numerical evaluation of the induced masses.

The softmin parameter \(\lambda\) introduces an additional approximation specific to the smoothed supporting-hyperellipsoid construction. As \(\lambda\to0\), the smoothed representation approaches the original nonsmooth supporting-hyperellipsoid representation, whereas larger values of \(\lambda\) yield smoother approximations at the expense of a larger smoothing error. Thus, for the implementation considered here, the approximation is influenced by the target discretization, the numerical evaluation of the reflector-induced masses, and the softmin smoothing, while the computational cost is primarily determined by \(K\), \(N\), and the number of supporting-hyperellipsoid iterations.

Several directions remain open for further investigation, including adaptive construction of the discrete target measure, more efficient evaluation of the reflector-induced masses, and alternative numerical methods for constructing the reflecting surface within the GOAS framework.

\subsection{Conclusions}

We have introduced Geometric Optics Approximation Sampling (GOAS), a transport-based sampling framework that constructs a deterministic map from a simple source measure to a prescribed target measure. As a concrete realization, we developed a supporting-hyperellipsoid construction that requires only a discrete approximation of the target measure and no gradient information of the target density. The formulation accommodates both density-based and sample-based target approximations. A softmin smoothing was further introduced to obtain a smooth approximate transport map.

We established the well-posedness and stability of the resulting push-forward measure and derived quantitative error estimates in the MMD metric, and established convergence of continuous statistical observables, including fixed-order moments.

Numerical experiments validated the theoretical results and demonstrated the applicability of GOAS to analytically tractable targets, non-Gaussian distributions, sample-based target approximations, and Bayesian inverse problems.




\begin{appendices}
The appendices includes mathematical formulation of near-field reflector problem, details on numerical experiments, additional numerical experiments, and all proof details.

\section{Details on numerical experiments}

\subsection{Spherical reflecting surface}\label{AppendixSph}
Let $n=2$, and then the $\Omega\subset P=\{x_3=h|h<0\}$ and $\varGamma\subset S^2_+$.
Assuming that the reflecting surface is a spherical sheet, i.e.,
\begin{align}\label{RefSph}
R(m)=rm,\quad m\in\varGamma,
\end{align}
where $r$ is the polar radius and a positive constant, then the direction of reflection is
$
y=-m.
$
Indeed, the unit normal direction $\upsilon$ of the reflecting surface is equal to the direction $m$ of ray emission.
Then from the \eqref{RefMap}, the reflection mapping is given by
\begin{align}\label{SphT}
z=T(m)=\frac{h}{w}m.
\end{align}
Hence from \eqref{PDE1} and \eqref{SphT}, we obtain
$
\pi(T(m))=\frac{I(m)w^3}{h^2}, \quad m\in\varGamma.
$
Then the density of target domain is given by
\begin{align}\label{LSph}
\pi(z)=\frac{-hI(-\hat{z})}{|z|^3},\quad z\in\Omega,
\end{align}
and the $\Omega$ is in the plane $P$ and depends on the $\varGamma$.
Indeed, if $\varGamma$ is assumed to be a spherical cap on the $S^2_+$ and its height is $h_c$,
then the $\Omega$ is a disk with radius $-h\sqrt{{1}/{(1-h_c)^2}-1}$.
And the normalisation constant of density \eqref{LSph} is given by $2\pi Ih_c$ if the $I$ is a uniform distribution.


\subsection{Strongly non-Gaussian distributions}\label{subsec:GOASvsMCMCs}
We compare our GOAS method with triangular transport map (TM) methods \cite{ricardo2024,parno2022} and several standard MCMCs: Metropolis-Hastings (MH) \citep{robert2004}, slice sampling \citep{neal2003}, Hamiltonian Monte Carlo (HMC) \citep{duane1987,neal2011}, and Metropolis-Adjusted Langevin Algorithm (MALA) \citep{roberts1998}.
The evaluation employs five challenging two-dimensional synthetic distributions--Funnel, Banana, Mixture of Gaussians (MoG), Ring, and Cosine \citep{wenliang2019,jaini2019}--which collectively exhibit diverse geometric structures and multimodality.

For GOAS, all target distributions are discretized using Hammersley point sets. For the MCMCs, we generate \(2\times10^5\) samples and retain \(5\times10^4\) samples after burn-in. For the TM method, the map is parameterized using a probabilists' Hermite polynomial basis of total degree \(p=5\), and \(2\times10^4\) samples are generated from the learned map.

\begin{figure}[htbp]
 \centering
 \includegraphics[width=0.7\textwidth,height=0.6\textwidth]{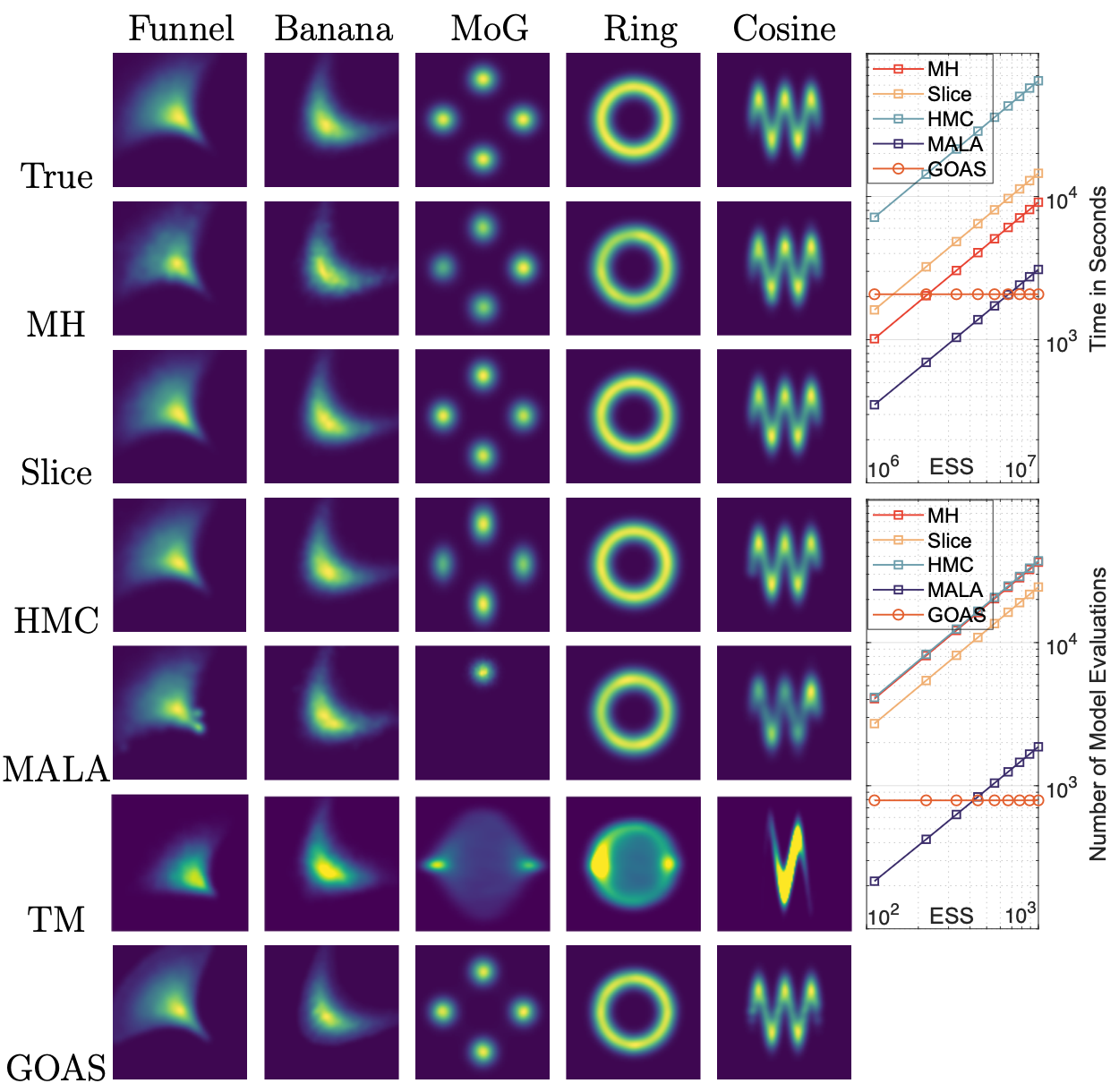}
 \caption{
Comparison of GOAS with a transport-map (TM) and MCMCs for non-Gaussian target distributions. The first five columns show the true densities and kernel density estimates obtained from MCMCs, TM and GOAS. The right panels show computational time and the number of density evaluations versus effective sample size (ESS) for the MoG distribution sampling.
}
\label{fig:GOASvsMCMCs}
\end{figure}

Figure~\ref{fig:GOASvsMCMCs} compares the true densities with kernel density estimates obtained from MH, slice sampling, HMC, MALA, TM, and GOAS. GOAS reproduces the main geometric and modal features of all five target distributions. The MCMCs also provide accurate approximations in several cases, although their performance varies with the geometry of the target. The polynomial TM performs well for the Funnel and Banana distributions, but the fixed-order map has greater difficulty representing the multimodal, ring-shaped, and oscillatory structures in the MoG, Ring, and Cosine examples. In contrast, GOAS captures these structures without directly parameterizing the transport map in a prescribed polynomial or neural-network class. Once the reflecting surface is constructed, the corresponding map is explicitly determined by the physical law of reflection. Moreover, the GOAS construction does not require derivatives of the target density, unlike gradient-based samplers such as HMC and MALA.

The right panels of Figure~\ref{fig:GOASvsMCMCs} compare computational time and the number of density evaluations as functions of the effective sample size (ESS) for the MoG example. For the MCMCs, both quantities increase with ESS. In contrast, the dominant computational cost of GOAS is incurred once in constructing the transport map, after which additional independent samples can be generated at comparatively small additional cost. Details of the ESS calculation are provided in Appendix~\ref{AppendixESSHD}.

\subsubsection{Evaluation of effective sample size}\label{AppendixESSHD}
Here we describe the calculation of effective sample size (ESS) \citep{liu2001monte}.
Let the $\tau_i$ be the integrated autocorrelation time of dimension $i$, and it is given by
$$
\tau_i=1+2\sum_{j=1}^{N_s}corr(\theta_1,\theta_{1+j})
$$
for dimension $i$ of samples $\{\theta_j\}_{j=1}^{N_s}$, where $corr(\cdot,\cdot)$ is the correlation coefficient.
Then we define the maximum integrated autocorrelation time over all dimension
$
\tau_{max}=\max_{i\in\{1,2,\dots,n\}}\tau_i.
$
The ESS is then computed by
$
ESS=\frac{N_s}{\tau_{max}}.
$

\subsection{Bayesian framework and numerical experiments}
\subsubsection{Bayesian framework}\label{Appendix:bayfra}
Let $\mathbf{X}$ be a separable Hilbert space, equipped with the Borel $\sigma$-algebra, and
 $\mathcal{G}:\mathbf{X}\to \mathbb{R}^n$ be a measurable function called the forward operator, which represents the connection between parameter and data in the mathematical model.
We wish to solve the inverse problems of finding the unknow model parameters $u$ in set $\mathbf{X}$ from measurement data $y\in\mathbb{R}^n$, which is usually generated by
\begin{align}\label{ForMod}
y=\mathcal{G}(u)+\eta,
\end{align}
where the noise $\eta$ is assumed to be a $n$-dimensional zero-mean Gaussian random variable with covariance matrix $\Sigma_{\eta}$.
In Bayesian Bayesian inference, the unknown model input and the measurement data are usually regarded as random variables \citep{stuart2010inverse}.
From \eqref{ForMod}, we define the negative log-likelihood
$$
\Phi_y(u):=\frac{1}{2}|\Sigma_{\eta}^{-\frac{1}{2}}(\mathcal{G}(u)-y)|^2.
$$
Combining the prior probability measure $\mu_0$ with density $\pi_0$ and Bayes theorem gives the posterior density up to a normalizing constant
\begin{align}\label{PosDen}
\pi(u) = \exp\bigl(-\Phi_y(u)\bigr)\pi_0(u).
\end{align}

\subsubsection{Locating acoustic source }\label{Appendix:detailsLAS}

The parameters are set as $\varsigma=(1,1,1)$, $k=1$, with $180$ directions on $S^1$.
The acoustic sources in \eqref{Hel} are located at $z_1=(1,4),z_2=(2,5),z_3=(3,6)$ and let $\vartheta=(1,2,3,4,5,6)$.
We employ a Gaussian prior $\mu_0$ with mean $(3,3,3,3,3,3)$ and identity covariance matrix.
Measurement data are obtained by $y=\mathcal{G}(\vartheta)+\eta$ where $\mathcal{G}$ represents the far field pattern \eqref{FarFie} and $\eta$ is the Gaussian noise with the standard deviation equal to noise level $5\%$ of the maximum norm of the forward model output.
We ran $10^5$ MCMC samples, discarding the first $3\times 10^4$ samples and using the remaining samples as the reference MCMC estimates for this example.

The target measure for GOAS is the corresponding Bayesian posterior measure. To construct a discrete approximation of the posterior, we employ an importance-sampling strategy. A pilot MCMC ($10^4$ iterations) run is first used only to estimate the posterior mean and covariance, denoted by \(\widehat m\) and \(\widehat C\), respectively. These estimates are used to define a Gaussian proposal
$
q(z)=\mathcal N(z;\widehat m,\widehat C).
$
We then generate \(K\) independent candidate points $z_1,\dots,z_K$ from $q$ and assign the unnormalized importance weights
\[
\widetilde w_i
=
\frac{\pi(z_i)}{q(z_i)},
\]
where \(\widetilde\pi\) denotes the unnormalized posterior density. After normalization,
\[
w_i
=
\frac{\widetilde w_i}
{\sum_{j=1}^K\widetilde w_j},
\]
the posterior is approximated by the weighted discrete measure
\[
\mu_t^K
=
\sum_{i=1}^K w_i\delta_{z_i}.
\]

\subsubsection{Initial-field reconstruction in a nonlinear advection--diffusion--reaction model}\label{Appendix:detailsADR}

We use the finite element method with Newton iteration to solve the equation \eqref{ADR}-\eqref{ADR3}.
We take $T=0.005,\; \Omega=[0,1]\times [0,1],\; c=1,\; \nu=(cos(t),sin(t)),\;m_d=0.01x_1$ and
$
f(x,t)=\exp |x-0.5|^2 / 0.9^2.
$
The number of time steps of discretization is $21$ and the $\Omega$ is discretised into a triangular mesh with $328$ elements and $185$ vertices.
Let $x=(x_1,x_2)$ and the basis of space $L^2(\Omega)$ is truncated as
$
\{cos(2\pi nx_1)cos(2\pi mx_2)+cos(2\pi nx_1)sin(2\pi mx_2)+sin(2\pi nx_1)cos(2\pi mx_2)+sin(2\pi nx_1)sin(2\pi mx_2)\}_{n,m=0}^{\mathbb{N},\mathbb{M}}
$
where $\mathbb{N},\mathbb{M}$ are integers.
In this example we inverse the coefficients of $m_0$ in this trigonometric basis.
The exact initial condition field in \eqref{ADR}is set
\begin{align*}
m_0(x)=&2+0.2cos(2\pi x_2)+0.3sin(2\pi x_2)+0.4cos(2\pi x_1)+0.5sin(2\pi x_1)\\
&+0.6cos(2\pi x_1)cos(2\pi x_2)+0.7cos(2\pi x_1)sin(2\pi x_2)\\
&+0.8sin(2\pi x_1)cos(2\pi x_2)+0.9sin(2\pi x_1)sin(2\pi x_2).
\end{align*}
We specify a Gaussian prior $\mu_0$ with $[1,0.15ones(1,8)]$ mean and $diag([1,0.1ones(1,8)]^2)$ covariance matrix.
The measurement data are obtained by $y=\mathcal{G}(m_0)+\eta$ where $\mathcal{G}$ is the forward model \eqref{ADR}-\eqref{ADR3} and $\eta$ represents the Gaussian noise with the standard deviation taken by noise level $1\%$ of the maximum norm of $u(x,T)$.
To avoid `inverse crimes', we generate measurement data by solving the forward problem on a finer grid.
We ran $10^5$ MCMC samples, discarding the first $5\times 10^4$ samples and using the remaining samples as the reference MCMC estimates for this example.

For this example, we construct the discrete posterior using a low-discrepancy point set. A pilot MCMC chain of length \(10^4\) is first used only to identify a bounded computational region
\[
\mathcal D
=
\prod_{\ell=1}^{9}[a_\ell,b_\ell]
\subset\mathbb R^9
\]
containing the dominant posterior mass. We then generate a Hammersley point set
$
\{z_j\}_{j=1}^{N_q}\subset\mathcal D
$
and evaluate the unnormalized posterior density \(\pi(z_j)\) at these points. To avoid retaining points with negligible posterior mass, we keep only those satisfying
\[
\frac{\pi(z_j)}
{\max_{1\le \ell\le N_q}\pi(z_\ell)}
\ge \epsilon.
\]
Let \(z_1,\ldots,z_K\) denote the retained points. Since the Hammersley points provide an approximately uniform discretization of \(\mathcal D\), their target weights are defined by
\[
w_i
=
\frac{\pi(z_i)}
{\displaystyle\sum_{j=1}^{K}\pi(z_j)},
\qquad
i=1,\ldots,K.
\]
The resulting discrete posterior measure is therefore
\[
\mu_t^K
=
\sum_{i=1}^{K}w_i\delta_{z_i},
\]
which is subsequently used to construct the GOAS reflector. Notice that this discretization requires only pointwise evaluations of the unnormalized posterior density; neither its normalizing constant nor its gradient is required.

\section{Mathematical formulation of near-field reflector problem}\label{sec:matref}
We derive the governing equations for a reflector surface that transforms a source distribution into a desired near-field target distribution. Assuming mass conservation, we have
\begin{align}\label{EneCon}
 \int_{E} I(m) d\sigma(m)=\int_{T(E)} \pi(z) \diff\mu(z)
\end{align}
for any open set $E\subset \varGamma$.
Note that
$$
\int_{T(E)} \pi(z) d\mu(z)=\int_{E} \pi\bigl(T(m)\bigr) \bigl|\det\bigl(J(T)\bigr)\bigr|\diff \sigma(m),
$$
where $\det\bigl(J(T)\bigr)$ is the Jacobian determinant of the map $T$.
Therefore, we obtain
\begin{align}\label{PDE1}
\bigl|\det\bigl(J(T(m))\bigr)\bigr|=\frac{I(m)}{\pi\bigl(T(m)\bigr)},\quad m\in \varGamma.
\end{align}
A necessary compatibility condition is
\begin{align}\label{TotEneCon}
 \int_{\varGamma} I(m) d\sigma(m)=\int_{\Omega} \pi(z) \diff\mu(z).
\end{align}

Let $(t^1,t^2,\dots,t^n)$ be a smooth parametrization of $S^n$.
Then $m=m(t^1,t^2,\dots,t^n)$.
Denote by $(e_{ij})=(\partial_i m \cdot \partial_j m)$ the matrix of coefficients of the first fundamental forms of $S^n$, where $\partial_i=\partial/\partial t^i$.
Put $(e^{ij})=(e_{ij})^{-1}$ and $\nabla=e^{ij}\partial_jm\partial_i$.
Then the unit normal vector of $R(m)=m\rho(m)$ is given by
\begin{align}\label{NorRef}
\upsilon=\frac{\nabla \rho -m\rho}{\sqrt{\rho^2+|\nabla \rho|^2}},
\end{align}
where $|\nabla \rho|^2=e^{ij}\partial_j\rho \partial_i\rho$.
Indeed, the unit normal vector \eqref{NorRef} can be obtain by $\partial_1R \times \partial_2R\times\dots\times\partial_nR$ where $\times$ denotes the outer product in space $\mathbb{R}^{n+1}$.
For simplicity in computing the Jacobian determinant of $T$, we choose an coordinate system on $S^n$ near the north pole.
Let $m=(m_1,m_2,\cdots,m_{n+1})\in \varGamma$ satisfying
\begin{align}\label{CooSys}
\begin{cases}
m_k(t^1,t^2,\dots,t^n) &=t^k,\quad k=1,2,\dots,n\\
m_{n+1}(t^1,t^2,\dots,t^n) & =w,
\end{cases}
\end{align}
and $|t|<1$, where $w:=\sqrt{1-|t|^2}$ and $t=(t^1,t^2,\dots,t^n)\in \Omega_{\varGamma} \subset \mathbb{R}^n$.
Therefore we also regard $\rho=\rho(t)$ as a function on $\Omega_{\varGamma}$ and $T=T(t)$ as a mapping on $\Omega_{\varGamma}$.

\begin{thm}\label{theoremPDE}
Let $\varGamma \subset S^n_+$ and $\Omega\subset P$, and let the density functions $I\in L^1(\varGamma)$ and $\pi\in L^1(\Omega)$ be given, satisfying the energy conservation \eqref{EneCon}.
For the polar radius $\rho$ of reflecting surface $R$ in the reflector shape design problem, define $u:=1/\rho$, and
$$
    a:=|Du|^2-(u-Du\cdot t)^2, \quad b:=|Du|^2+u^2-(Du\cdot t)^2,
$$
where $Du=(\partial_1u,\partial_2u,\dots,\partial_nu)$ is the gradient of $u$, and
$$
\mathcal{M}:=1+\frac{t\otimes t}{1-|t|^2},\quad c:=\frac{h}{w},
$$
where $t\otimes t=(t^it^j)$ is an $n\times n$ matrix.
Then the $u$ is governed by the equation
\begin{align}\label{AM}
 \det\biggl(D^2u+\frac{ca}{2(1-cu)}\mathcal{M}\biggr)=\frac{|a^{n+1}|}{2^n(1-cu)^nb}\cdot \frac{I(t)}{w\pi\circ T(t)},\quad t\in\Omega_{\varGamma},
\end{align}
where $D^2u=(\partial_i\partial_j u)$ is the Hessian matrix of $u$.
\end{thm}
The equation \eqref{AM} is a fully nonlinear elliptic partial differential equation of Monge-Amp\`{e}re type, and the proof follows a similar approach as in \citep{karakhanyan2010}. 



\begin{rem}
\begin{enumerate}[(i)]
\item Substituting $u=1/\rho$ into equation \eqref{AM}, we can obtain the equation
\begin{align}\label{PDErho}
 \det\biggl(-D^2\rho+\frac{2}{\rho}D\rho\otimes D\rho+\frac{c\tilde{a}}{2\rho(\rho -c)}\mathcal{M}\biggr)=\frac{|\tilde{a}^{n+1}|}{2^n\rho^n(\rho-c)^n\tilde{b}}\cdot \frac{I(t)}{w\pi\circ T(t)}
\end{align}
for $\rho$, where $\tilde{a}=|D\rho|^2-(\rho+D\rho\cdot t)^2$ and $\tilde{b}=|D\rho|^2+\rho^2-(D\rho\cdot t)^2$.
The boundary condition for this equation is
\begin{align}\label{eq:PDEboundary}
T(\varGamma)=\Omega.
\end{align}

\item If one set $\Omega\subset\{x_{n+1}=0\}$, i.e., $h=0$, then the \eqref{AM} simplifies to
\begin{align}\label{StaAM}
 \det(D^2u)=\frac{|a^{n+1}|}{2^nb}\cdot \frac{I(t)}{w\pi\circ T(t)},
\end{align}
which is a standard Monge-Amp\`{e}re equation.
\end{enumerate}
\end{rem}

\begin{proof}
Let the reflection map $T=(Ts,T_{n+1})$, and denote $T_{i,j}=\partial_jT_i, i=1,2,\dots,n+1, \text{and}\, j=1,2,\dots,n$.
Let $\eta=(0,0,\dots,1)$ denote the unit normal vector of $\Omega$.
We get the Jacobian determinant of $T$
\begin{align}\label{detJT}
\det(J(T))&=\frac{(\partial_1T(m)\times\partial_2T(m)\times\cdots \times\partial_nT(m))\cdot \eta}{\sqrt{det(e_{ij})}}\notag\\
&=w
\begin{vmatrix}
T_{1,1} & T_{1,2} & \ldots & T_{1,n}&\eta_1\\
T_{2,1} & T_{2,2} & \ldots & T_{2,n}&\eta_2\\
\vdots & \vdots & \ddots & \vdots & \vdots\\
T_{n,1} & T_{n,2} & \ldots & T_{n,n}&\eta_n\\
T_{n+1,1} & T_{n+1,2} & \ldots & T_{n+1,n}&\eta_{n+1}\\
\end{vmatrix}\notag\\
&=w
\begin{vmatrix}
T_{1,1} & T_{1,2} & \ldots & T_{1,n}\\
T_{2,1} & T_{2,2} & \ldots & T_{2,n}\\
\vdots & \vdots & \ddots & \vdots \\
T_{n,1} & T_{n,2} & \ldots & T_{n,n}\\
\end{vmatrix}\notag\\
&=:w\det(DTs).
\end{align}

From the \eqref{NorRef} and using the orthonormal coordinate system \eqref{CooSys}, we obtain the  unit normal vector of $R$
\begin{align}\label{NorRef1}
\upsilon=\frac{(D\rho,0) -m(t\cdot D\rho+\rho)}{\sqrt{\rho^2+|D\rho|^2-(t\cdot D\rho)^2}},
\end{align}
Then
$$
m\cdot \upsilon=-\frac{\rho}{\sqrt{\rho^2+|D\rho|^2-(D\rho\cdot m)^2}}.
$$
Hence, the reflection direction is
\begin{align}\label{RefDirP}
y&=m-2(m\cdot \upsilon)\upsilon \notag\\
&=m+\frac{2\rho((D\rho,0)-m(\rho+D\rho\cdot m))}{\rho^2+|D\rho|^2-(D\rho\cdot m)^2}\notag\\
&=m\frac{|D\rho|^2-(\rho+D\rho\cdot m)^2}{\rho^2+|D\rho|^2-(D\rho\cdot m)^2}+\frac{2\rho(D\rho,0)}{\rho^2+|D\rho|^2-(D\rho\cdot m)^2}.
\end{align}

Since $\Omega\subset P$ and \eqref{RefMap}, we have
\begin{align*}
 m_{n+1}\rho+y_{n+1}l=h
\end{align*}
and from \eqref{RefDirP},
\begin{align}\label{yn1}
y_{n+1}=m_{n+1}\frac{|D\rho|^2-(\rho+D\rho\cdot m)^2}{\rho^2+|D\rho|^2-(D\rho\cdot m)^2},
\end{align}
therefore,
\begin{align}\label{dp}
l=\biggl(\frac{h}{w}-\rho\biggr)\frac{\rho^2+|D\rho|^2-(D\rho\cdot m)^2}{|D\rho|^2-(\rho+D\rho\cdot m)^2}.
\end{align}
By \eqref{RefMap},\eqref{RefDirP} and \eqref{dp}, we get
\begin{align}\label{ComT}
T=\frac{-2\rho^2(D\rho,0)}{|D\rho|^2-(\rho+D\rho\cdot m)^2}+\biggl(m-\frac{2\rho(D\rho,0)}{|D\rho|^2-(\rho+D\rho\cdot m)^2}\biggr)\frac{h}{w}.
\end{align}
Then
\begin{align}\label{Ts}
Ts=\frac{2Du}{|Du|^2-(u-Du\cdot t)^2}+\biggl(t-\frac{2uDu}{|Du|^2-(u-Du\cdot t)^2}\biggr)\frac{h}{w}.
\end{align}
Consider a general mapping $Q$ from $\Omega_{\varGamma}\times\mathbb{R}\times\mathbb{R}^n$ into $\mathbb{R}^n$, and denote points in $\Omega_{\varGamma}\times\mathbb{R}\times\mathbb{R}^n$ by $(x,r,p)$, we see from \eqref{Ts}
\begin{align}\label{YW}
Q(x,r,p)&=\frac{2p}{|p|^2-(r-p\cdot x)^2}+\biggl(x-\frac{2rp}{|p|^2-(r-p\cdot x)^2}\biggr)\frac{h}{w}\notag\\
&:=Y+W,
\end{align}
and then
\begin{align}\label{DTs}
DTs&=DQ\notag\\
&=Q_pD^2u+Q_x+Q_r\otimes Du\notag\\
&=(Y_p+W_p)\bigl[D^2u+(Y_p+W_p)^{-1}(Y_x+Y_r\otimes Du+W_x+W_r\otimes Du)\bigr]\notag\\
&:=(Y_p+W_p)\bigl[D^2u+A(\cdot,r,p)\bigr].
\end{align}
By \eqref{PDE1}, \eqref{detJT} and \eqref{DTs}, we then have
\begin{align}\label{PDE2}
\Bigl|\det\bigl[D^2u+A(\cdot,r,p)\bigr]\Bigr|=\frac{1}{\bigl|\det[Y_p+W_p]\bigr|}\cdot\frac{I}{wL\circ T}.
\end{align}
For $Y=Y(x,r,p)$ in \eqref{YW}, it is easy to compute that
\begin{align}\label{Yxr}
Y_x+Y_r\otimes Du=0,
\end{align}
and
\begin{align}\label{Yp}
Y_p=\frac{2}{|p|^2-(r-p\cdot t)^2}\Biggl[Id-\frac{2p\otimes \bigl[p+(r-p\cdot t)t\bigr]}{|p|^2-(r-p\cdot t)^2}\Biggr].
\end{align}
Similarly,
\begin{align}\label{Wxr}
W_x+W_r\otimes Du=cId+\frac{ct\otimes t}{1-|t|^2}-\frac{2c}{|p|^2-(r-p\cdot t)^2}(\frac{rp\otimes t}{1-|t|^2}+p\otimes p),
\end{align}
and
\begin{align}\label{Wp}
W_p=-crY_p.
\end{align}
Hence,
\begin{align}\label{YWp}
Y_p+W_p=\frac{2(1-cr)}{a}\Biggl[Id-\frac{2p\otimes \bigl[p+(r-p\cdot t)t\bigr]}{a}\Biggr].
\end{align}
Using the formula $\det[Id+\alpha\otimes\beta]=1+\alpha\cdot\beta$ and  $(Id+\alpha\otimes\beta)^{-1}=Id-\alpha\otimes\beta/(1+\alpha\cdot\beta)$ for any vector $\alpha,\beta \in\mathbb{R}^n$, we have
\begin{align}\label{detYWp}
\det[Y_p+W_p]=-\frac{2^n(1-cr)^nb}{a^{n+1}},
\end{align}
and
\begin{align}\label{invYWp}
(Y_p+W_p)^{-1}=\frac{a}{2(1-cr)}\Biggl[Id-\frac{2p\otimes \bigl[p+(r-p\cdot t)t\bigr]}{b}\Biggr].
\end{align}
Hence, by \eqref{Yxr}-\eqref{YWp} and \eqref{invYWp}, we get
\begin{align}\label{Ares}
A(\cdot,r,p)=\frac{ca}{2(1-cr)}(1+\frac{t\otimes t}{1-|t|^2}).
\end{align}
Combining \eqref{PDE2}, \eqref{detYWp} and \eqref{Ares}, we have then obtain the equation \eqref{AM}.
\end{proof}

Monge-Amp\`{e}re type partial differential equations commonly appear in optimal transport problems, where the goal is to find a mass-preserving transport map that minimizes a given cost functional \citep{villani2009,villani2021}. However, the near-field reflector problem differs fundamentally from classical optimal transport. While it induces a transportation mapping, the associated cost function depends nonlinearly on the potential, and its weak solutions do not optimize the corresponding cost functional \citep{karakhanyan2010,graf2012}.

\section{Proof Details}
This supplementary section contains the detailed proofs of the theorems and lemmas stated in the main text. The purpose is to provide a complete and rigorous justification of the mathematical framework and the theoretical results underlying the proposed method.

\subsection{Proof in Section \ref{sec:geobay}}

\begin{defn}[Convex reflecting surface]\label{def:ConRefSur}
A reflecting surface $R_\rho$ is called convex with respect to $\Omega$ if every $m\in\varGamma$ admits a supporting hyperellipsoid. 
If $R_\rho$ is composed of $K$ hyperellipsoidal patches,
\begin{align}\label{Rpolyhedron}
    R_\rho
    =
    \bigcup_{i=1}^K (E_{z_i}\cap R_\rho),
    \qquad z_i\in\Omega,
\end{align}
we call it a polyhedral reflector.
\end{defn}

For a convex reflecting surface, the polar radius can be represented as the lower envelope of supporting hyperellipsoids.

\begin{prop}\label{rhoTconvexRef}
Let $R_\rho$ be a convex reflecting surface. 
Then its polar radius satisfies
\begin{align}\label{ConRefRho}
    \rho(x)
    =
    \inf_{z\in\Omega}
    \frac{d(z)}
    {1-\varepsilon(d(z))(\hat z\cdot x)},
    \qquad x\in\varGamma .
\end{align}
\end{prop}

\begin{proof}
From the Definition \ref{def:ConRefSur}, for any $x\in\varGamma$, there is a supporting hyperellipsoid at $x\rho(x)$ whose polar radius is denoted by $\rho_z(x)$.
By the Definition \ref{def:SE}, we have
$$
\rho(x)=\rho_z(x),\quad \rho(x')\leq \rho_z(x'),\quad \forall x'\in\varGamma,
$$
and then \eqref{ConRefRho} is obtained.
\end{proof}

\subsubsection*{Proof of Theorem \ref{StaRef}}
\begin{proof}
From the measure $\lim_{k\to\infty}\mu_t^k=\mu_t$, there exists a monotone convergent sequence $\{\Omega_k\}_{k=1}^{\infty}$, i.e., $\Omega_1\supseteq \Omega_2\dots\supseteq \Omega_k\supseteq \Omega_{k+1}\dots$, and $\lim_{k\to\infty}\Omega_k=\cap_{k=1}^{\infty}\Omega_k=\Omega$ such that
$$
\mu_t^k=\int_{\Omega_k}\pi(x)\diff\mu(x),\quad \mu_t=\int_\Omega \pi(x)\diff\mu(x).
$$
By \eqref{ConRefRho}, for each the convex reflector $R_{\rho_k}$ corresponding to measure $\mu_t^k$, the polar radius is given by
$$
\rho_k(x)=\inf_{z\in \Omega_k}\frac{d(z)}{1-\varepsilon(d(z))(\hat{z}\cdot x)},\quad x\in \varGamma.
$$
Hence,
\begin{align*}
\lim_{k\to\infty}\rho_k(x)&=\lim_{k\to\infty}\inf_{z\in \Omega_k}\frac{d(z)}{1-\varepsilon(d(z))(\hat{z}\cdot x)}\\
&=\inf_{z\in \Omega}\frac{d(z)}{1-\varepsilon(d(z))(\hat{z}\cdot x)}\\
&=\rho(x),
\end{align*}
for any fixed $x\in\varGamma$.
That is, we get the pointwise convergence, i.e., $\rho_k(x)\to\rho(x),\forall x\in\varGamma$.
Noting the monotonicity of sequence $\{\Omega_k\}_{k=1}^{\infty}$ and infimum functional, and combining this with the continuity of $\rho_k$ and $\rho$, we get that $\rho_k$ is monotone with respect to $k$.
Then, by the Dini theorem, $\rho_k$ converges uniformly to $\rho$, namely
$$
\|\rho_k-\rho\|_{C(\varGamma)}\to 0,\quad k\to\infty.
$$
By the relation between the Hausdorff metric and the Lipschitz norm \citep{Groemer1994,graf2012}, we have
$$
d_{H}(R_{\rho_k},R_{\rho})\leq\|\rho_k-\rho\|_{C(\varGamma)},
$$
which completes the proof.
\end{proof}

\subsubsection*{Proof of Theorem \ref{GOAMstab}}
\begin{proof}
For each $\mu_t^k$, let $R_k$ be the corresponding reflecting surface satisfying
$$
\mu_{R_k}(w)=\mu_t^k(w),\quad w\subseteq\Omega.
$$
The result follows directly from Theorem \ref{StaRef}.
\end{proof}

\subsection{Proof in Section \ref{sec:algorithm}}
\subsubsection*{Proof of Theorem \ref{thm:rhoSmo}}

\begin{proof}
From \eqref{ConRefRho} we obtain \eqref{eq:rhodis}.
If there exits $m\geq 1$ minimal points with $f_{z_{i_k}}(x;d_{i_k})=\tilde{\rho}(x),k\in I^{\ast}$, we have $\tau_j\geq\tau\geq 0$ for any $j\not\in I^{\ast}$.
Then,
\begin{align*}
\tilde{\rho}_{\lambda}(x)&=-\lambda\log\biggl[\sum_{i=1}^{K}e^{-\frac{f_{z_i}(x;d_i)}{\lambda}}\biggr]\\
&=-\lambda\log\biggl[\sum_{i\in I^{\ast}}e^{-\frac{\tilde{\rho}(x)}{\lambda}}+\sum_{i\not\in I^{\ast}}e^{-\frac{f_{z_i}(x;d_i)}{\lambda}}\biggr]\\
&=\tilde{\rho}(x)-\lambda\log\biggl(m+\sum_{i\not\in I^{\ast}}e^{-\frac{\tau_i}{\lambda}}\biggr).
\end{align*}
Thus,
\begin{align*}
|\tilde{\rho}_{\lambda}-\tilde{\rho}|&=\biggl|\lambda\log\biggl(m+\sum_{i\not\in I^{\ast}}e^{-\frac{\tau_i}{\lambda}}\biggr)\biggr|\\
&\leq \lambda\log(m+(K-m)e^{-\frac{\tau}{\lambda}})\\
&=\lambda\biggl[\log(m)+\log(1+\frac{K-m}{m}e^{-\frac{\tau}{\lambda}})\biggr]\\
&\leq \lambda\log(m)+\frac{\lambda(K-m)}{m}e^{-\frac{\tau}{\lambda}},
\end{align*}
which establishes \eqref{ineq:rhoappro}.

If $m=1$, let $I^{\ast}=\{i^{\ast}\}$.
For $i\neq i^{\ast}$,
$$
e^{-\frac{f_{z_i}(x;d_i)}{\lambda}}=e^{-\frac{\tau_i}{\lambda}}e^{-\frac{f_{z_{i^{\ast}}}(x;d_{i^{\ast}})}{\lambda}}\leq e^{-\frac{\tau}{\lambda}}\sum_{i=1}^{K}e^{-\frac{f_{z_i}(x;d_i)}{\lambda}}.
$$
Since $f$ is smooth,
\begin{align*}
|D\tilde{\rho}_{\lambda}-D\tilde{\rho}|&=\biggl|\frac{\sum_{i=1}^{K}e^{-\frac{f_{z_i}(x;d_i)}{\lambda}}Df_{z_i}(x;d_i)}{\sum_{i=1}^{K}e^{-\frac{f_{z_i}(x;d_i)}{\lambda}}}-Df_{z_{i^{\ast}}}(x;d_{i^{\ast}})\biggr|\\
&=\biggl|\frac{\sum_{i\neq i^{\ast}}e^{-\frac{f_{z_i}(x;d_i)}{\lambda}}(Df_{z_i}(x;d_i)-Df_{z_{i^{\ast}}}(x;d_{i^{\ast}}))}{\sum_{i=1}^{K}e^{-\frac{f_{z_i}(x;d_i)}{\lambda}}}\biggr|\\
&\leq \frac{\sum_{i\neq i^{\ast}}e^{-\frac{\tau}{\lambda}}\sum_{i=1}^{K}e^{-\frac{f_{z_i}(x;d_i)}{\lambda}}|Df_{z_i}(x;d_i)-Df_{z_{i^{\ast}}}(x;d_{i^{\ast}})|}{\sum_{i=1}^{K}e^{-\frac{f_{z_i}(x;d_i)}{\lambda}}}\\
&\leq e^{-\frac{\tau}{\lambda}}\sum_{i\neq i^{\ast}}|Df_{z_i}(x;d_i)-Df_{z_{i^{\ast}}}(x;d_{i^{\ast}})|\\
&\leq C(K-1)e^{-\frac{\tau}{\lambda}},
\end{align*}
where positive parameter $C$ depends on $f$.
\end{proof}

\subsubsection*{Proof of Proposition \ref{prop:refmap}}

\begin{proof}
From \eqref{NorRef} and using the orthonormal coordinate system \eqref{CooSys}, the unit normal of $R_K$ is given by
\begin{align}\label{eq:unitnormalapp}
\upsilon=\frac{D\tilde{\rho}_{\lambda}-x(\tilde{\rho}_{\lambda}+D\tilde{\rho}_{\lambda}\cdot x)}{\sqrt{\tilde{\rho}_{\lambda}^2+|D\tilde{\rho}_{\lambda}|^2-(D\tilde{\rho}_{\lambda}\cdot x)^2}}.
\end{align}
Then
$$
x\cdot \upsilon=-\frac{\tilde{\rho}_{\lambda}}{\sqrt{\tilde{\rho}_{\lambda}^2+|D\tilde{\rho}_{\lambda}|^2-(D\tilde{\rho}_{\lambda}\cdot x)^2}}.
$$
The reflection direction becomes
\begin{align}\label{eq:RefDirDis}
y=x-2(x\cdot \upsilon)\upsilon
=x+\frac{2\tilde{\rho}_{\lambda}(D\tilde{\rho}_{\lambda}-x(\tilde{\rho}_{\lambda}+D\tilde{\rho}_{\lambda}\cdot x))}{\sqrt{\tilde{\rho}_{\lambda}^2+|D\tilde{\rho}_{\lambda}|^2-(D\tilde{\rho}_{\lambda}\cdot x)^2}}.
\end{align}
Since $\Omega\subset P$ and by \eqref{RefMap}, we have
$
x_{n+1}\tilde{\rho}_{\lambda}+y_{n+1}l=z_{n+1},
$
and from \eqref{eq:RefDirDis}
\begin{align}\label{}
y_{n+1}=x_{n+1}\frac{|D\tilde{\rho}_{\lambda}|^2-(\tilde{\rho}_{\lambda}+D\tilde{\rho}_{\lambda}\cdot x)^2}{\tilde{\rho}_{\lambda}^2+|D\tilde{\rho}_{\lambda}|^2-(D\tilde{\rho}_{\lambda}\cdot x)^2}.
\end{align}
Thus,
\begin{align}\label{eq:lDis}
l=\biggl(\frac{z_{n+1}}{x_{n+1}}-\tilde{\rho}_{\lambda}\biggr)\frac{\tilde{\rho}_{\lambda}^2+|D\tilde{\rho}_{\lambda}|^2-(D\tilde{\rho}_{\lambda}\cdot x)^2}{|D\tilde{\rho}_{\lambda}|^2-(\tilde{\rho}_{\lambda}+D\tilde{\rho}_{\lambda}\cdot x)^2}.
\end{align}
Combining \eqref{RefMap}, \eqref{eq:RefDirDis} and \eqref{eq:lDis}, we then obtain the equation \eqref{eq:DisT}.
\end{proof}

\subsection{Proof in Section \ref{sec:errorEstimation}}
\begin{thm}
\label{thm:MMD_IS}
Let $\nu_t$ be a measure on $\Omega$, and let $\nu$ be another measure such that $\nu_t$ is absolutely continuous with respect to $\nu$. Let $X_1, \dots, X_N$ be i.i.d. (i.e., independent and identically distributed) random samples from $\nu$. Define the importance weight $w(x) := \frac{\diff\nu_t}{\diff\nu}(x)$, and the empirical measure
\[
\nu_t^N := \frac{1}{N} \sum_{i=1}^N w(X_i) \delta_{X_i}.
\]
Assume $\mathbb{E}_\nu[w(X)^2] < \infty$ and assumption \ref{assump:kernel_bounded} holds. Then
\[
\mathbb{E}_\nu [d_k(\nu_t^N, \nu_t)^2] = \frac{1}{N} \left( \mathbb{E}_\nu[w(X)^2 k(X, X)] - \mathbb{E}_{\nu_t}[k(X, X')] \right),
\]
where $X, X'$ are independent with law $\nu_t$. In particular,
\[
\mathbb{E}_\nu [d_k(\nu_t^N, \nu_t)] \le \frac{\kappa}{\sqrt{N}} \sqrt{\mathbb{E}_\nu[w(X)^2]}.
\]
\end{thm}

\begin{proof}
Let $m_{\nu_t} := \int k(x, \cdot) \diff\nu_t(x)$ and $m_{\nu_t^N} := \int k(x, \cdot) \diff\nu_t^N(x)$ be the kernel mean embeddings in the RKHS $\mathcal{H}_k$. By definition of $\nu_t^N$, we have
\[
m_{\mu^N} = \frac{1}{N} \sum_{i=1}^N w(X_i) k(X_i, \cdot),
\]
which is well-defined in $\mathcal{H}_k$, because of assumptions $\mathbb{E}_\nu[w(X)^2] < \infty$ and \ref{assump:kernel_bounded}.
The MMD satisfies $d_k(\nu_t^N, \nu_t) = \|m_{\nu_t^N} - m_{\nu_t}\|_{\mathcal{H}_k}$. Define the random variables $\xi_i \in \mathcal{H}_k$ as
\[
\xi_i := w(X_i) k(X_i, \cdot) - m_{\nu_t}.
\]
Using the change of measure $\diff\nu_t = w \diff\nu$, we observe that
\[
\mathbb{E}_\nu[w(X_i) k(X_i, \cdot)] = \int w(x) k(x, \cdot) d\nu(x) = \int k(x, \cdot) \diff\nu_t(x) = m_{\nu_t},
\]
which implies $\mathbb{E}_\nu[\xi_i] = 0$. Since $X_i$ are i.i.d., we have
\begin{align*}
\mathbb{E}_\nu [d_k(\nu_t^N, \nu_t)^2] &= \mathbb{E}_\nu \left\| \frac{1}{N} \sum_{i=1}^N \xi_i \right\|_{\mathcal{H}_k}^2 \\
&= \frac{1}{N^2} \sum_{i,j=1}^N \mathbb{E}_\nu \langle \xi_i, \xi_j \rangle_{\mathcal{H}_k}.
\end{align*}
For $i \neq j$, the independence and zero-mean property give $\mathbb{E}_\nu \langle \xi_i, \xi_j \rangle = \langle \mathbb{E}_\nu \xi_i, \mathbb{E}_\nu \xi_j \rangle = 0$. Thus, only the diagonal terms remain:
\[
\mathbb{E}_\nu [d_k(\nu_t^N, \nu_t)^2] = \frac{1}{N} \mathbb{E}_\nu \|\xi_1\|_{\mathcal{H}_k}^2.
\]
Expanding the squared norm, we obtain
\[
\mathbb{E}_\nu \|\xi_1\|^2_{\mathcal{H}_k} = \mathbb{E}_\nu \|w(X_1) k(X_1, \cdot)\|^2_{\mathcal{H}_k} - 2 \mathbb{E}_\nu \langle w(X_1) k(X_1, \cdot), m_{\nu_t} \rangle + \|m_{\nu_t}\|^2_{\mathcal{H}_k}.
\]
By the reproducing property, $\|w(X_1) k(X_1, \cdot)\|^2_{\mathcal{H}_k} = w(X_1)^2 k(X_1, X_1)$. Furthermore,
\[
\mathbb{E}_\nu \langle w(X_1) k(X_1, \cdot), m_{\nu_t} \rangle = \langle \mathbb{E}_\nu [w(X_1) k(X_1, \cdot)], m_{\nu_t} \rangle = \langle m_{\nu_t}, m_{\nu_t} \rangle = \|m_{\nu_t}\|^2_{\mathcal{H}_k}.
\]
Substituting these back, and noting that $\|m_{\nu_t}\|^2_{\mathcal{H}_k} = \mathbb{E}_{\nu_t}[k(X, X')]$, we have
\[
\mathbb{E}_\nu \|\xi_1\|^2_{\mathcal{H}_k} = \mathbb{E}_\nu [w(X)^2 k(X, X)] - \|m_{\nu_t}\|^2_{\mathcal{H}_k} = \mathbb{E}_\nu [w(X)^2 k(X, X)] - \mathbb{E}_{\nu_t}[k(X, X')],
\]
which proves the first identity. 

Finally, since $\mathbb{E}_{\nu_t}[k(X, X')] \ge 0$ and assumption \ref{assump:kernel_bounded} holds, we have $\mathbb{E}_\nu [d_k({\nu_t}^N, {\nu_t})^2] \le \frac{\kappa^2}{N} \mathbb{E}_\nu [w(X)^2]$. Applying Jensen's inequality, $\mathbb{E}_\nu d_k \le \sqrt{\mathbb{E}_\nu d_k^2}$, yields the desired bound:
\[
\mathbb{E}_\nu [d_k(\nu_t^N, \nu_t)] \le \frac{\kappa}{\sqrt{N}} \sqrt{\mathbb{E}_\nu [w(X)^2]}.
\]
\end{proof}

In the special case $\nu=\nu_t$ (i.e., $w\equiv1$), this result reduces to the standard Monte Carlo convergence of empirical kernel mean embeddings \citep{gretton2012,muandet2017}.


\begin{lem}\label{lem:MMD_pullback_mean_embedding}
Let $T:\varGamma \to \Omega$ be measurable, and let $k:\Omega\times\Omega\to\mathbb{R}$ be a kernel with RKHS $\mathcal H_k$.
Define the pullback kernel
\[
k_T(x,y) := k(T(x),T(y))
\]
with RKHS $\mathcal H_{k_T}$.
Then, for any measures $\mu_1,\mu_2$ on $\varGamma$,
\[
d_k(T_\sharp\mu_1,T_\sharp\mu_2)
=
d_{k_T}(\mu_1,\mu_2).
\]
\end{lem}

\begin{proof}
For $i=1,2$, the mean embedding of $T_\sharp\mu_i$ in $\mathcal H_k$ is
\[
m_{T_\sharp\mu_i}
=
\int_\Omega k(z, \cdot)\, d(T_\sharp\mu_i)(z).
\]
By change of variables,
\[
m_{T_\sharp\mu_i}
=
\int_\Gamma k(T(x), \cdot)\, d\mu_i(x).
\]
Since the mean embedding of $\mu_i$ in $\mathcal H_{k_T}$ is given by
\[
m_{\mu_i}^{(k_T)}
=
\int_\Gamma k(T(x), \cdot)\, d\mu_i(x),
\]
we have
\[
m_{\mu_i}^{(k_T)}=m_{T_\sharp\mu_i}.
\]
Thus
\[
d_k(T_\sharp\mu_1,T_\sharp\mu_2)
=
\|m_{T_\sharp\mu_1}-m_{T_\sharp\mu_2}\|_{\mathcal H_k}
=
\|m_{\mu_1}^{(k_T)}-m_{\mu_2}^{(k_T)}\|_{\mathcal H_{k_T}}
=
d_{k_T}(\mu_1,\mu_2).
\]
\end{proof}

\begin{lem}\label{lem:MMD_map}
Let $\nu$ be a measure on $\varGamma$, and let $T_1,T_2:\varGamma\to\Omega$ be measurable maps. Assume $k$ satisfies assumption \ref{assump:kernel_lipschitz}. 
Then
\begin{equation}\label{eq:MMD_map}
d_{k}\bigl((T_1)_\sharp\nu,(T_2)_\sharp\nu\bigr)
\;\le\;
L_f\,\|T_1-T_2\|_{L^1(\nu)},
\end{equation}
where $L_f$ is a positive constant.
\end{lem}

\begin{proof}
By the change of variables,
\begin{align*}
d_{k}((T_1)_\sharp\nu,(T_2)_\sharp\nu)
&=\sup_{\|f\|_{\mathcal{H}_k} \leq 1} 
\left| 
\int_{\Omega} f \, \diff(T_1)_\sharp\nu - \int_{\Omega} f \, \diff(T_2)_\sharp\nu 
\right| \\
&=\sup_{\|f\|_{\mathcal{H}_k}\le1}
\left|
\int_{\varGamma} \bigl(f(T_1(x))-f(T_2(x))\bigr)\,\diff\nu(x)
\right|.
\end{align*}
For any $f\in\mathcal{H}_k$ with $\|f\|_{\mathcal{H}_k}\le1$, the reproducing property and assumption \ref{assump:kernel_lipschitz} imply
\[
|f(z_1)-f(z_2)| \le L_f |z_1-z_2|, \qquad \forall z_1,z_2\in\Omega.
\]
Hence
\[
|f(T_1(x))-f(T_2(x))|
\le L_f|T_1(x)-T_2(x)|.
\]
Integrating and taking the supremum yields \eqref{eq:MMD_map}.
\end{proof}

\begin{lem}\label{lem:Tdis}
For $\varGamma\subset S_+^n$ and $\Omega\subset P=\{z_{n+1}=h,h<0\}$, let $\widetilde{T}_{\lambda}:\varGamma\to\Omega$ be the approximate reflection map under $\widetilde{\rho}_{\lambda}$ in \eqref{eq:rhodisSmooth}, and let $\widetilde{T}:\varGamma\to\Omega$ be the reflection map under $\widetilde{\rho}$ in \eqref{eq:rhodis}.
Then
$$
|\widetilde{T}_{\lambda}-\widetilde{T}|\leq C_1\lambda e^{-\frac{1}{\lambda}}+C_2e^{-\frac{1}{\lambda}}.
$$
where $C_1,C_2$ are positive constants.
\end{lem}
\begin{proof}
From \eqref{RefMap}, we have
\begin{align*}
\widetilde{T}_{\lambda}-\widetilde{T}&=x\tilde{\rho}_{\lambda}+l_{\lambda}(x)y_{\lambda}(x)-(x\tilde{\rho}+l(x)y(x))\\
&=x(\tilde{\rho}_{\lambda}-\tilde{\rho})+y_{\lambda}(l_{\lambda}-l)+l(y_{\lambda}-y)\\
&=I_1+I_2+I_3.
\end{align*}
By \eqref{ineq:rhoappro} and $x\in S^n$,
$$
|I_1|\leq |x||\tilde{\rho}_{\lambda}-\tilde{\rho}|\leq C_1\lambda e^{-\frac{1}{\lambda}}.
$$

From \eqref{RefDir},
\begin{align*}
|I_3|&\leq|l||y_{\lambda}-y|\\
&\leq C_2|(x\cdot \upsilon)\upsilon-(x\cdot \upsilon_{\lambda})\upsilon_{\lambda}|\\
&\leq C_2|x\cdot(\upsilon-\upsilon_{\lambda})\upsilon+(x\cdot \upsilon_{\lambda})(\upsilon-\upsilon_{\lambda})|.
\end{align*}
Using \eqref{eq:unitnormalapp},\eqref{ineq:rhoappro}, \eqref{ineq:Drhoappro}, and $\upsilon,\upsilon_{\lambda}\in S^n$, we obtain
\begin{align*}
|\upsilon_{\lambda}-\upsilon|&\leq C_3|(D\tilde{\rho}_{\lambda}-D\tilde{\rho})-x((\tilde{\rho}_{\lambda}-\tilde{\rho})+(D\tilde{\rho}_{\lambda}-D\tilde{\rho})\cdot x)|\\
&\leq C_3(2|D\tilde{\rho}_{\lambda}-D\tilde{\rho}|+|\tilde{\rho}_{\lambda}-\tilde{\rho}|)\\
&\leq C_3C_4e^{-\frac{1}{\lambda}}+C_3\lambda e^{-\frac{1}{\lambda}},
\end{align*}
where $C_3,C_4>0$.
Thus,
$$
|I_3|\leq 2C_2|\upsilon_{\lambda}-\upsilon|\leq 2C_2(C_3C_4e^{-\frac{1}{\lambda}}+C_3\lambda e^{-\frac{1}{\lambda}}).
$$

Since $\Omega\subset P$ and $y_{\lambda}\in S^n$,
\begin{align*}
|I_2|&\leq|l_{\lambda}-l|\\
&=\biggl|\frac{z_{n+1}-x_{n+1}\tilde{\rho}_{\lambda}}{y_{\lambda,n+1}}-\frac{z_{n+1}-x_{n+1}\tilde{\rho}}{y_{n+1}}\biggr|\\
&=\biggl|\frac{x_{n+1}(\tilde{\rho}-\tilde{\rho}_{\lambda})}{y_{\lambda,n+1}}+\frac{(z_{n+1}-x_{n+1}\tilde{\rho})(y_{n+1}-y_{\lambda,n+1})}{y_{\lambda,n+1}y_{n+1}}\biggr|\\
&\leq C_5\lambda e^{-\frac{1}{\lambda}}+C_6|y_{n+1}-y_{\lambda,n+1}|\\
&\leq C_5\lambda e^{-\frac{1}{\lambda}}+ 2C_6(C_3C_4e^{-\frac{1}{\lambda}}+C_3\lambda e^{-\frac{1}{\lambda}}).
\end{align*}

Combining these estimates, we conclude
$$
|\widetilde{T}_{\lambda}-\widetilde{T}|\leq|I_1|+|I_2|+|I_3|,
$$
yielding the desired bound.
\end{proof}

\begin{thm}\label{thm:errorBoundDis}
Let $\widetilde{T}_{\lambda}$ and $\widetilde{T}$ be defined as in Lemma~\ref{lem:Tdis}. The discrete form of source measure and target measure are given by
\[
\mu_s^N := \frac{1}{N} \sum_{i=1}^N\delta_{x_i},
\]
where random samples $x_i$ from source measure, and $\mu_t^K$ in \eqref{eq:DisTarDen}, respectively.
Then
\[
d_{k}\bigl((\widetilde{T}_{\lambda})_\sharp\mu_s^N,
\mu_t^K\bigr)
\;\le\;
C_1\bigl(\lambda e^{-1/\lambda}+e^{-1/\lambda}\bigr)+C_2\epsilon,
\]
where $C_1,C_2$ are positive constant, and $\epsilon$ is any in advance error bound in \eqref{ErrBonDis}.
\end{thm}

\begin{proof}
By Lemma~\ref{lem:MMD_map} and Lemma~\ref{lem:Tdis},
\[
d_{k}\bigl((\widetilde{T}_{\lambda})_\sharp\mu_s^N,
\widetilde{T}_\sharp\mu_s^N\bigr)
\le
L_f\int_{\varGamma}|\widetilde{T}_{\lambda}-\widetilde{T}|\diff\mu_s^N
\le
C_1\bigl(\lambda e^{-1/\lambda}+e^{-1/\lambda}\bigr).
\]
Under \eqref{ErrBonDis} and \eqref{ApprG}, we have
\begin{align*}
d_{k}(\widetilde{T}_{\sharp}\mu_s^N,\mu_t^K)&\le\left| \int_{\Omega} f \, \diff\widetilde{T}_{\sharp}\mu_s^N - \int_{\Omega} f \, \diff\mu_t^K \right|\\
&\le \int_{\Omega} \|f\|_{\mathcal{H}_k} |\widetilde G-w|\diff\mu\\
&\le C_2M\epsilon
\end{align*}
By the triangle inequality, we have
\begin{align*}
d_{k}((\widetilde{T}_{\lambda})_{\sharp}\mu_s^N,\mu_t^K)&\leq d_{k}((\widetilde{T}_{\lambda})_{\sharp}\mu_s^N,\widetilde{T}_{\sharp}\mu_s^N)+d_{k}(\widetilde{T}_{\sharp}\mu_s^N,\mu_t^K)\\
&\leq C_1(\lambda e^{-\frac{1}{\lambda}}+e^{-\frac{1}{\lambda}})^{\frac{1}{2}} +C_2(M\epsilon),
\end{align*}
which completes the proof.
\end{proof}

\subsubsection*{Proof of Theorem \ref{thm:errorBound}}
\begin{proof}
From the triangle inequality and Lemma~\ref{lem:MMD_pullback_mean_embedding}:
\begin{align*}
d_k\big((\widetilde T_\lambda)_\sharp \mu_s, \mu_t\big)
&\le 
d_k\big((\widetilde T_\lambda)_\sharp \mu_s, (\widetilde T_\lambda)_\sharp \mu_s^N\big)
+ d_k\big((\widetilde T_\lambda)_\sharp \mu_s^N, \mu_t\big) \\
&\le 
d_{k_{\widetilde T_\lambda}}(\mu_s,\mu_s^N)
+ d_k\big((\widetilde T_\lambda)_\sharp \mu_s^N,\mu_t^K\big)
+ d_k(\mu_t^K,\mu_t).
\end{align*}
By Assumption~\ref{assump:kernel_bounded}, there exists $\kappa>0$ such that
\[
\sup_{z\in\Omega} k(z,z)\le \kappa^2.
\]
Since $\widetilde T_\lambda(\Gamma)\subset \Omega$, for every $x\in\Gamma$,
\[
k_{\widetilde T_\lambda}(x,x)
=
k(\widetilde T_\lambda(x),\widetilde T_\lambda(x))
\le \kappa^2.
\]
Hence the pullback kernel $k_{\widetilde T_\lambda}$ is bounded.
Applying Theorem~\ref{thm:MMD_IS} with $\nu_s=\mu_s$, we have
$w_s=\frac{d\mu_s}{d\nu_s}\equiv1$, and thus
$\mathbb E_{\mu_s}[w_s^2]=1$. Therefore,
\[
\mathbb E_{\mu_s}
\big[
d_{k_{\widetilde T_\lambda}}(\mu_s,\mu_s^N)
\big]
\le
\frac{\kappa}{\sqrt N}.
\]
The middle term
\(
d_k\big((\widetilde T_\lambda)_\sharp \mu_s^N,\mu_t^K\big)
\)
is controlled by Theorem~\ref{thm:errorBoundDis}.
Combining the three estimates and taking expectations gives the desired bound.
\end{proof}

\subsubsection*{Proof of Corollary \ref{cor:distarget}}
\begin{proof}
Applying Theorem~\ref{thm:MMD_IS} to $\mu_t$ and $\mu_t^K$, with importance weight
$w_t=\frac{d\mu_t}{d\nu_t}$ and assumption
$\mathbb E_{\nu_t}[w_t(X)^2]<C_3^2$, yields
\[
\mathbb E_{\nu_t}
\big[
d_k(\mu_t^K,\mu_t)
\big]
\le
C_3\,\frac{\kappa}{\sqrt K}.
\]
\end{proof}

\subsubsection*{Proof of Theorem \ref{thm:moment_error}}
\begin{proof}
Since $k$ is universal, $\mathcal H_k$ is dense in $C(\Omega)$ with respect to the supremum norm. Hence, for any $\delta>0$ and $p\in C(\Omega)$, there exists $f_\delta\in\mathcal H_k$ such that
\[
\|p-f_\delta\|_\infty\le\delta.
\]
Therefore,
\begin{align*}
&
\left|
\int_\Omega p\,\diff\mu_{\rm GOAS}
-
\int_\Omega p\,\diff\mu_t
\right|
\\
&\le
\left|
\int_\Omega (p-f_\delta)\,\diff\mu_{\rm GOAS}
\right|
+
\left|
\int_\Omega (p-f_\delta)\,\diff\mu_t
\right|
+
\left|
\int_\Omega f_\delta\,\diff\mu_{\rm GOAS}
-
\int_\Omega f_\delta\,\diff\mu_t
\right|.
\end{align*}
Since both measures have total mass $M$,
\[
\left|
\int_\Omega (p-f_\delta)\,\diff\mu_{\rm GOAS}
\right|
\le M\delta,
\qquad
\left|
\int_\Omega (p-f_\delta)\,\diff\mu_t
\right|
\le M\delta.
\]
By the definition of MMD,
\[
\left|
\int_\Omega f_\delta\,\diff\mu_{\rm GOAS}
-
\int_\Omega f_\delta\,\diff\mu_t
\right|
\le
\|f_\delta\|_{\mathcal H_k}
d_k(\mu_{\rm GOAS},\mu_t).
\]
Combining the above estimates gives
\[
\left|
\int_\Omega p\,\diff\mu_{\rm GOAS}
-
\int_\Omega p\,\diff\mu_t
\right|
\le
2M\delta
+
\|f_\delta\|_{\mathcal H_k}
d_k(\mu_{\rm GOAS},\mu_t).
\]

Now suppose that $d_k(\mu_{\rm GOAS},\mu_t)\to0$. For any $\varepsilon>0$, choose $\delta>0$ such that
$2M\delta<\varepsilon/2$. For this fixed $\delta$, the function $f_\delta$ is fixed, and hence $d_k(\mu_{\rm GOAS},\mu_t)$ can be chosen sufficiently small so that
\[
\|f_\delta\|_{\mathcal H_k}
d_k(\mu_{\rm GOAS},\mu_t)
<\frac{\varepsilon}{2}.
\]
Thus,
\[
\left|
\int_\Omega p\,\diff\mu_{\rm GOAS}
-
\int_\Omega p\,\diff\mu_t
\right|
<\varepsilon,
\]
which proves the result.
\end{proof}

\end{appendices}

\bibliographystyle{siamplain}
 \bibliography{bibliography.bib}

\begin{thebibliography}{10}

\bibitem{aronszajn1950}
{\sc N.~Aronszajn}, {\em Theory of reproducing kernels}, Transactions of the
  American Mathematical Society, 68 (1950), pp.~337--404.

\bibitem{ricardo2024}
{\sc R.~Baptista, Y.~Marzouk, and O.~Zahm}, {\em On the representation and
  learning of monotone triangular transport maps}, Foundations of Computational
  Mathematics, 24 (2024), pp.~2063--2108.

\bibitem{berlinet2011}
{\sc A.~Berlinet and C.~Thomas-Agnan}, {\em Reproducing Kernel Hilbert Spaces
  in Probability and Statistics}, Kluwer Academic Publishers, 2004.

\bibitem{steve2011}
{\sc S.~Brooks, A.~Gelman, G.~L. Jones, and X.-L. Meng}, eds., {\em Handbook of
  Markov Chain Monte Carlo}, Chapman and Hall/CRC, 2011.

\bibitem{duane1987}
{\sc S.~Duane, A.~D. Kennedy, B.~J. Pendleton, and D.~Roweth}, {\em Hybrid
  monte carlo}, Physics Letters B, 195 (1987), pp.~216--222.

\bibitem{el2011}
{\sc A.~El~Badia and T.~Nara}, {\em An inverse source problem for helmholtz's
  equation from the cauchy data with a single wave number}, Inverse Problems,
  27 (2011), p.~105001.

\bibitem{el2012}
{\sc T.~A. El~Moselhy and Y.~M. Marzouk}, {\em Bayesian inference with optimal
  maps}, Journal of Computational Physics, 231 (2012), pp.~7815--7850.

\bibitem{eller2009}
{\sc M.~Eller and N.~P. Valdivia}, {\em Acoustic source identification using
  multiple frequency information}, Inverse Problems, 25 (2009), p.~115005.

\bibitem{fournier2010}
{\sc F.~R. Fournier, W.~J. Cassarly, and J.~P. Rolland}, {\em Fast freeform
  reflector generation using source-target maps}, Optics Express, 18 (2010),
  pp.~5295--5304.

\bibitem{gelman2013}
{\sc A.~Gelman, J.~B. Carlin, H.~S. Stern, D.~B. Dunson, A.~Vehtari, and D.~B.
  Rubin}, {\em Bayesian Data Analysis}, CRC Press, 3~ed., 2013.

\bibitem{ghattas2021}
{\sc O.~Ghattas and K.~Willcox}, {\em Learning physics-based models from data:
  Perspectives from inverse problems and model reduction}, Acta Numerica, 30
  (2021), pp.~445--554.

\bibitem{glassner1989}
{\sc A.~S. Glassner}, ed., {\em An Introduction to Ray Tracing}, Morgan
  Kaufmann, 1989.

\bibitem{graf2012}
{\sc T.~Graf and V.~I. Oliker}, {\em An optimal mass transport approach to the
  near-field reflector problem in optical design}, Inverse Problems, 28 (2012),
  p.~025001.

\bibitem{gretton2012}
{\sc A.~Gretton, K.~M. Borgwardt, M.~J. Rasch, B.~Sch{\"o}lkopf, and A.~Smola},
  {\em A kernel two-sample test}, Journal of Machine Learning Research, 13
  (2012), pp.~723--773.

\bibitem{Groemer1994}
{\sc H.~Groemer}, {\em Stability results for convex bodies and related
  spherical integral transformations}, Advances in Mathematics, 109 (1994),
  pp.~45--74.

\bibitem{jaini2019}
{\sc P.~Jaini, K.~A. Selby, and Y.~Yu}, {\em Sum-of-squares polynomial flow},
  in Proceedings of the 36th International Conference on Machine Learning,
  vol.~97 of Proceedings of Machine Learning Research, PMLR, 2019,
  pp.~3009--3018.

\bibitem{karakhanyan2010}
{\sc A.~L. Karakhanyan and X.-J. Wang}, {\em On the reflector shape design},
  Journal of Differential Geometry, 84 (2010), pp.~561--610.

\bibitem{kaufman2009finding}
{\sc L.~Kaufman and P.~J. Rousseeuw}, {\em Finding groups in data: an
  introduction to cluster analysis}, John Wiley \& Sons, 2009.

\bibitem{kobyzev2020}
{\sc I.~Kobyzev, S.~J.~D. Prince, and M.~A. Brubaker}, {\em Normalizing flows:
  An introduction and review of current methods}, IEEE Transactions on Pattern
  Analysis and Machine Intelligence, 43 (2020), pp.~3964--3979.

\bibitem{kochengin1997}
{\sc S.~A. Kochengin and V.~I. Oliker}, {\em Determination of reflector
  surfaces from near-field scattering data}, Inverse Problems, 13 (1997),
  pp.~363--373.

\bibitem{kochengin1998}
{\sc S.~A. Kochengin and V.~I. Oliker}, {\em Determination of reflector
  surfaces from near-field scattering data ii. numerical solution}, Numerische
  Mathematik, 79 (1998), pp.~553--568.

\bibitem{liu2001monte}
{\sc J.~S. Liu}, {\em Monte Carlo Strategies in Scientific Computing},
  Springer, 2001.

\bibitem{liu2021}
{\sc Y.~Liu, Y.~Guo, and J.~Sun}, {\em A deterministic-statistical approach to
  reconstruct moving sources using sparse partial data}, Inverse Problems, 37
  (2021), p.~065005.

\bibitem{marzouk2016}
{\sc Y.~Marzouk, T.~Moselhy, M.~Parno, and A.~Spantini}, {\em Sampling via
  measure transport: An introduction}, in Handbook of Uncertainty
  Quantification, Springer, 2016, pp.~1--41.

\bibitem{muandet2017}
{\sc K.~Muandet, K.~Fukumizu, B.~Sriperumbudur, and B.~Sch{\"o}lkopf}, {\em
  Kernel mean embedding of distributions: A review and beyond}, Foundations and
  Trends in Machine Learning, 10 (2017), pp.~1--141.

\bibitem{neal2003}
{\sc R.~M. Neal}, {\em Slice sampling}, The Annals of Statistics, 31 (2003),
  pp.~705--767.

\bibitem{neal2011}
{\sc R.~M. Neal}, {\em Mcmc using hamiltonian dynamics}, in Handbook of Markov
  Chain Monte Carlo, Chapman and Hall/CRC, 2011, pp.~113--162.

\bibitem{papamakarios2021}
{\sc G.~Papamakarios, E.~Nalisnick, D.~J. Rezende, S.~Mohamed, and
  B.~Lakshminarayanan}, {\em Normalizing flows for probabilistic modeling and
  inference}, Journal of Machine Learning Research, 22 (2021), pp.~1--64.

\bibitem{park2009simple}
{\sc H.-S. Park and C.-H. Jun}, {\em A simple and fast algorithm for k-medoids
  clustering}, Expert systems with applications, 36 (2009), pp.~3336--3341.

\bibitem{parno2022}
{\sc M.~Parno, P.-B. Rubio, D.~Sharp, M.~Brennan, R.~Baptista, H.~Bonart, and
  Y.~Marzouk}, {\em {MParT}: Monotone parameterization toolkit}, Journal of
  Open Source Software, 7 (2022), p.~4843.

\bibitem{rezende2015}
{\sc D.~J. Rezende and S.~Mohamed}, {\em Variational inference with normalizing
  flows}, in Proceedings of the 32nd International Conference on Machine
  Learning, vol.~37 of Proceedings of Machine Learning Research, PMLR, 2015,
  pp.~1530--1538.

\bibitem{robert2004}
{\sc C.~P. Robert and G.~Casella}, {\em Monte Carlo Statistical Methods},
  Springer, 2~ed., 2004.

\bibitem{roberts1998}
{\sc G.~O. Roberts and J.~S. Rosenthal}, {\em Optimal scaling of discrete
  approximations to langevin diffusions}, Journal of the Royal Statistical
  Society: Series B, 60 (1998), pp.~255--268.

\bibitem{Schmitzer2019}
{\sc B.~Schmitzer}, {\em Stabilized sparse scaling algorithms for entropy
  regularized transport problems}, SIAM Journal on Scientific Computing, 41
  (2019), pp.~A1443--A1481.

\bibitem{peter2008}
{\sc P.~Shirley and R.~K. Morley}, {\em Realistic Ray Tracing}, AK Peters,
  2008.

\bibitem{stuart2010inverse}
{\sc A.~M. Stuart}, {\em Inverse problems: A bayesian perspective}, Acta
  Numerica, 19 (2010), pp.~451--559.

\bibitem{tierney1994}
{\sc L.~Tierney}, {\em Markov chains for exploring posterior distributions},
  The Annals of Statistics, 22 (1994), pp.~1701--1762.

\bibitem{villani2021}
{\sc C.~Villani}, {\em Topics in Optimal Transportation}, vol.~58 of Graduate
  Studies in Mathematics, American Mathematical Society, 2003.

\bibitem{villani2009}
{\sc C.~Villani}, {\em Optimal Transport: Old and New}, vol.~338 of Grundlehren
  der mathematischen Wissenschaften, Springer, 2009.

\bibitem{wenliang2019}
{\sc L.~K. Wenliang, D.~J. Sutherland, H.~Strathmann, and A.~Gretton}, {\em
  Learning deep kernels for exponential family densities}, in Proceedings of
  the 36th International Conference on Machine Learning, vol.~97 of Proceedings
  of Machine Learning Research, PMLR, 2019, pp.~6737--6746.

\end{thebibliography}

\end{document}